\documentclass[11pt,reqno,a4paper]{amsart}
\usepackage[
    a4paper,
    top=2.5cm,
    bottom=2.5cm,
    left=2.5cm,
    right=2.5cm
]{geometry}
\usepackage{indentfirst, amssymb,amsmath,amsthm, mathrsfs,cite,graphicx,float}

\usepackage{tikz}
\usepackage{float}
\newtheorem*{theoA}{Theorem A}
\newtheorem*{theoB}{Theorem B}
\newtheorem*{theoC}{Theorem C}

\newtheorem{thm}{Theorem}[section]
\newtheorem{lem}{Lemma}[section]

\newtheorem{cor}{Corollary}[section]

\newtheorem{rem}{Remark}[section]

\newtheorem{defi}{Definition}[section]

\newtheorem{open problem}{Open problem}[section]
\newtheorem{ex}{Example}[section]

\newcommand{\be}{\begin{equation}}
\newcommand{\ee}{\end{equation}}
\newcommand{\bs}{\begin{small}}
\newcommand{\es}{\end{small}}
\newcommand{\beas}{\begin{eqnarray*}}
\newcommand{\eeas}{\end{eqnarray*}}
\newcommand{\bea}{\begin{eqnarray}}
\newcommand{\eea}{\end{eqnarray}}
\renewcommand{\epsilon}{\varepsilon}
\numberwithin{equation}{section}
\newcommand{\norm}[1]{\left\lVert#1\right\rVert}

\begin{document}
\title[Multidimensional mixed Bohr radii]{Multidimensional mixed Bohr radii for Banach space valued holomorphic functions}
\author[V. Allu and R. Biswas]{Vasudevarao Allu and Raju Biswas}
\date{}
\address{Vasudevarao Allu, Department of Mathematics, School of Basic Sciences, Indian Institute of Technology Bhubaneswar, Bhubaneswar-752050, Odisha, India.}
\email{avrao@iitbbs.ac.in}
\address{Raju Biswas, Department of Mathematics, Raiganj University, Raiganj, West Bengal-733134, India.}
\email{rajubiswasjanu02@gmail.com}
\maketitle
\let\thefootnote\relax
\begin{abstract} 
In this paper, we introduce and study a generalized notion of multidimensional mixed Bohr radii for holomorphic functions defined on the unit ball of the finite-dimensional Minkowski 
sequence space with values in complex Banach spaces. We first prove that the multidimensional mixed Bohr radius is positive for every non-null bounded linear operator whose norm 
is sufficiently small relative to the scaling parameter. We then determine the exact asymptotic behavior of the multidimensional mixed Bohr radius for infinite-dimensional Banach 
spaces of finite cotype and we obtain analogous bounds for spaces without finite cotype. Finally, for concave Banach lattices and summing operators, we prove lower bounds that 
involve a logarithmic factor. These results connect classical scalar Bohr phenomena with operator-valued and Banach space valued settings, and highlight how the geometry of the 
spaces and the mapping properties of the operator influence coefficient estimates and their dependence on 
dimension. Several examples illustrate the consistency of our results with the classical Bohr radius and the sharp Defant--Frerick asymptotics.

\end{abstract}

\section{Introduction}
 In the early years of the twentieth century, Harald Bohr \cite{B1914} made a significant contribution to the development of a general theory of Dirichlet series of the form:
\bea\label{b1}
	D(s):=\sum_{n\geq 1} a_n n^{-s},
\eea
where $a_n\in \Bbb{C}$ and $s=\sigma +it$ is a complex variable. Every Dirichlet series of the form \eqref{b1} can be viewed as a sequence of functions 
$\left(\sum_{n=1}^k a_n n^{-s}\right)_{k\in\Bbb{N}}$ and it converges on half-planes of the form $\{s\in \Bbb{C}: \mbox{Re}(s)>\sigma, \,\sigma\in\Bbb{R}\}$. It is well-known 
that the width of the largest possible strip on which the Dirichlet series \eqref{b1} converges but does not converge absolutely or uniformly, respectively, is $1$. 
Bohr was interested in determining the largest width of the strips in the plane where a Dirichlet series may converge uniformly without converging absolutely. This problem is 
well-known as \textit{Bohr's absolute convergence problem} in the theory of Dirichlet series. To tackle this problem, Bohr discovered a connection between a Dirichlet series and a power series of infinitely many variables. 
Using the fundamental theorem of arithmetic, Bohr considered the prime decomposition of $n=p_1^{\alpha_1}\cdots p_k^{\alpha_k}$ for each $n\in \Bbb{N}$ and defined $z=(p^{-s}_{1}, \ldots, p^{-s}_{k})$, where $p_1<p_2<\cdots<p_k$ are prime numbers. Thus,
\beas
	D(s)=\sum_{n\geq 1}a_n(p^{-s}_{1})^{\alpha_1}\cdots (p^{-s}_{k})^{\alpha_k}=\sum a_n z^{\alpha_1}_{1}\cdots z^{\alpha_k}_k.
\eeas
This correspondence is known as the \textit{Bohr transform}. It connects Dirichlet series and power series, enabling the use of complex analysis techniques to address problems involving Dirichlet series.\\[2mm]
\indent The Bohr phenomenon has its origins in a question about power series. Bohr was studying a power series in one complex variable and sought to determine the values of 
$z\in\Bbb{D}$ for which the inequality $\sum_{n=0}^{\infty}|a_n z^n| \leq \sup_{z\in \Bbb{D}}\left|\sum_{n=0}^{\infty}a_n z^n\right|$ holds. Since this inequality is not 
valid throughout the unit disk $\Bbb{D}$, a natural question arises: what is the largest subdomain of $\Bbb{D}$ on which the inequality holds for every power series that is 
holomorphic and bounded in $\Bbb{D}$? This question led to the development of the Bohr phenomenon.
In 1914, Bohr \cite{B1914} proved the following result, which is well-known as the {\it Bohr's inequality}.
\begin{theoA}\cite{B1914}
	Let $f(z)= \sum_{n=0}^{\infty}a_n z^n$ be a holomorphic function in the unit disk $\Bbb{D}$ such that $\sup_{z\in \Bbb{D}}|f(z)|<\infty$. Then, we have
	\beas
		\sum_{n=0}^{\infty}|a_n|r^n \leq \sup_{z\in \Bbb{D}}|f(z)|\quad\mbox{for $|z|=r \leq 1/3$}.
	\eeas
	 The constant $1/3$, known as the Bohr radius, cannot be replaced by a larger constant.
\end{theoA}
\indent
Bohr's original result gave the radius $1/6$, but this was soon sharpened to the optimal value $1/3$ by M. Riesz, I. Schur, and N. Wiener, each working independently. Bohr's paper \cite{B1914} presents both his own proof and Wiener's proof. Two decades later, Bohnenblust and Hille \cite{BH1931} resolved Bohr's absolute convergence problem by proving that 
the maximum possible width of such a strip is $1/2$. In 1995, Dixon \cite{D1995} provided an application of Bohr's inequality in Banach algebras, which motivated many 
mathematicians to work in this direction. For an overview of some aspects of the work on Bohr's inequality, we refer to a survey by Abu-Muhanna {\it et al.} \cite{AAP2017}. Subsequently, several versions of the Bohr inequality in the planar case have been extensively studied
(see \cite{KP2019, LLP2021, ABM2026, LP2021,KPW2025, DLP2025} and the references therein). \\[2mm]
\indent In 1997, Boas and Khavinson \cite{BK1997} introduced  and studied the $n$-dimensional Bohr radius $K_n$ for holomorphic functions on the $n$-dimensional polydisk 
$\Bbb{D}^n$. After that the field became very active, several researchers like Aizenberg and his coauthor \cite{A2000,AAD2000,AAD2001}, Balasubramanian {\it et al.} \cite{BCQ2006}, Allu {\it et al.} \cite{AHP2026}, Boas \cite{B2000}, Bombieri and Bourgain \cite{BB2004}, Galicer {\it et al.} \cite{GMM2020}, Hamada {\it et al.} \cite{HHK2009}, Paulsen {\it et al.} \cite{PPS2002}, Popescu \cite{P2019}, Kumar and Manna \cite{KM2025}, Kumar and Ponnusamy \cite{KP2026} have paid much attention to this result to use it in different contexts of mathematics and generalized it.

 
\section{Preliminaries}
We use standard terminology and notions from classical Banach space theory (see \cite{LT1977,LT1979}). Let $\Bbb{D}=\{z\in\Bbb{C} :|z|<1\}$ be the open unit disk in the
complex plane $\Bbb{C}$. Throughout this paper, it has been assumed that every Banach space $X$ is complex and $X^*$ is its topological dual. We write $B_{X}$ for the open unit ball in $X$. We denote $\ell^n_p$ for the Banach space of all $n$-tuples $z=(z_1, \ldots, z_n)\in \Bbb{C}^n$ endowed with the norm 
\beas
	\norm{z}_p:= \left(\sum_{i=1}^{n}|z_i|^p\right)^{\frac{1}{p}}\quad \mbox{for} \quad 1\leq p<\infty
\eeas
and $\norm{z}_{\infty}:= \max_{1\leq i\leq n}|z_i|$ for $p=\infty$. For $1\leq p\leq \infty$, the conjugate exponent $p'$ for $p$ is defined by $1/p+1/p'=1$. For $x\in X$, the absolute value of $x$ is defined by $|x|=x\vee (-x)$. A domain $\Omega \subseteq \mathbb{C}^n$ is 
said to be balanced if for every $z \in \Omega $ and every $\lambda \in \overline{\mathbb{D}}$, we have $\lambda z \in \Omega$. A Banach lattice is a Banach space $X$ which is a vector lattice with $|x|\leq |y|$ implies $\norm{x}\leq \norm{y}$ for all $x,y\in X$. A Banach lattice is said to be $r$-concave, where $1\leq r<\infty$, if there is a constant $L>0$ such that 
\bea\label{b3}
	\left(\sum_{j=1}^{n}\norm{x_j}^r\right)^{1/r}\leq L \norm{\left(\sum_{j=1}^{n}|x_j|^r\right)^{1/r}}
\eea 
for every finite family of elements $x_1,\ldots, x_n\in X$. The best such constant satisfying \eqref{b3} is denoted by $M_r(X)$. For $2\leq r\leq \infty$, a Banach space $X$ is said to have cotype $r$, if there is a constant $C>0$ such that for every choice of finitely many vectors $x_1, \ldots, x_n\in X$, we have
\bea\label{b4}
	\left(\sum_{j=1}^{n}\norm{x_j}^r\right)^{1/r} \leq C \left(\int_{0}^{1}\norm{\sum_{j=1}^{n}r_j(t)x_j}^2\, dt\right)^{1/2},
\eea
where $r_j$ is the $j$-th Rademacher function on $[0,1]$. Further, the best such constant $C$ satisfying \eqref{b4} is denoted by $C_r(X)$. We write 
\beas
	\text{Cot}(X):= \inf \left\{2\leq r\leq \infty: X \,\, \mbox{has cotype}\, r\right\}.
\eeas
Note that every Banach space has cotype $\infty$ and for $\text{Cot}(X)=\infty$, we denote $1/\text{Cot}(X)=0$. It is to be noted that the notions of concavity and cotype are 
closely related concepts in the context of Banach lattices. In fact, an $r$-concave Banach lattice with $r\geq 2$ is of cotype $r$. Conversely, each Banach lattice of cotype $2$ is 
$2$-concave. Furthermore, a Banach lattice of cotype $r>2$ is $t$-concave for all $t>r$.\\[2mm]
\indent Given $m,n\in \Bbb{N}$, we consider the following sets of indices:
\beas
	\mathcal{M}(m,n)&=&\left\{\mathbf{j}=(j_1,\ldots,j_m): 1\leq j_1,\ldots,j_m\leq n\right\}=\left\{1,2,\ldots,n\right\}^m,\\[2mm]
	\mathcal{J}(m,n)&=&\left\{\mathbf{j}\in \mathcal{M}(m,n): 1\leq j_1\leq \cdots\leq j_m\leq n\right\}\\[2mm]\text{and}\quad
	\Lambda(m,n)&=&\left\{\alpha=(\alpha_1,\ldots,\alpha_n)\in \Bbb{N}^n_{0}: |\alpha|:=\alpha_1+\cdots+\alpha_n=m\right\}.
\eeas
There is a one-to-one correspondence between the multi-index sets $\Lambda(m,n)$ and $\mathcal{J}(m,n)$. For each $\alpha=(\alpha_1,\ldots,\alpha_n)\in \Lambda(m,n)$, the associated index $\mathbf{j}_{\alpha}\in \mathcal{J}(m,n)$ is given by $$\mathbf{j}_{\alpha}=\left(1,\stackrel{\alpha_1}{\ldots},1, 2,\stackrel{\alpha_{2}}{\ldots}2, \ldots,n,\stackrel{\alpha_n}{\ldots},n\right).$$
Conversely, for each $\mathbf{j}:=(j_1,\ldots,j_m)\in \mathcal{J}(m,n)$ the associated index $\alpha:=(\alpha_1,\ldots,\alpha_n)\in \Lambda(m,n)$ is given by $\alpha_r=\left|\{k:j_k=r\}\right|$. 
For some fixed $\mathbf{j}:=(j_1,\ldots,j_m)\in \mathcal{J}(m,n)$, we define
 \beas [\mathbf{j}]=\left\{\mathbf{i}:=(i_1,\ldots,i_m)\in\Bbb{N}^m: (i_{\sigma(1)},\ldots,i_{\sigma(m)})=\mathbf{j}\quad\text{ for some permutation}\quad \sigma \in S_m \right\}\eeas
 and $|\mathbf{j}|=\text{card}[\mathbf{j}]$.
 We observe that
\begin{equation*}
	\mathcal{M}(m,n)=\cup_{\mathbf{j}\in \mathcal{J}(m,n)}[\mathbf{j}].
\end{equation*} 
Moreover, for each $\alpha\in \Lambda(m,n)$, we have
\beas
	\left|\mathbf{j}_{\alpha}\right|=\frac{m!}{\alpha!},\quad\text{where}\quad \mathbf{j}_{\alpha}=\left(1,\stackrel{\alpha_1}{\ldots},1, 2,\stackrel{\alpha_{2}}{\ldots}2, \ldots,n,\stackrel{\alpha_n}{\ldots},n\right)\quad\text{and}\quad \alpha!=\alpha_1!\cdots\alpha_n!.
\eeas
We consider the following notation from \cite{BDS2019}. Given a subset $\mathcal{J}\subset \mathcal{J}(m,n)$, we denote
\beas
	\mathcal{J}^*=\left\{\mathbf{j}\in \mathcal{J}(m-1,n): \,\mbox{there is}\,\,k\geq 1, (\mathbf{j},k)\in \mathcal{J}\right\}.
\eeas
For $1\leq p\leq \infty$, let $\mathcal{P}(^m\ell^n_p,X)$ be the Banach space of all $m$-homogeneous polynomial $P:\ell^n_p\rightarrow X$ equipped with the norm 
\beas
	\norm{P}_{\mathcal{P}(^m\ell^n_p,X)}:=\sup_{z\in B_{\ell^n_{p}}}\norm{P(z)}_{X}.
\eeas
Using the multi-index notations, a polynomial $P\in \mathcal{P}(^m\ell^n_p,X)$ can be written as
\beas
	P(z)=\sum_{\alpha\in\Lambda(m,n)}c_{\alpha}z^{\alpha},
\eeas
where $c_{\alpha}\in X$ and $z^\alpha=z_1^{\alpha_1} \cdot z_2^{\alpha_2}\cdots z_n^{\alpha_n}$ (where $0^0$ is represented as $1$).
There is another way to write a polynomial $P\in \mathcal{P}(^m\ell^n_p,X)$ as 
\beas
	P(z)=\sum_{\mathbf{j}\in \mathcal{J}(m,n)} a_{\mathbf{j}}z_{\mathbf{j}},
\eeas
where $a_{\mathbf{j}}\in X$, $z_{\mathbf{j}}=z_{j_1} z_{j_2}\cdots z_{j_m}$. Note that $c_{\alpha}=a_{\mathbf{j}}$ with
$\mathbf{j}=\left(1,\stackrel{\alpha_1}{\ldots},1, 2,\stackrel{\alpha_{2}}{\ldots}2, \ldots,n,\stackrel{\alpha_n}{\ldots},n\right)$. If $f:B_{\ell^n_p}\rightarrow X$ is an $X$-valued holomorphic function, then we write $f(z)=\sum_{\alpha\in \Bbb{N}^n_{0}}c_{\alpha}(f)z^{\alpha}$ as the monomial series expansion of $f$, where $c_{\alpha}(f)=(\partial^{\alpha}f(0))/\alpha!$ is the $\alpha$-th coefficient of the monomial series expansion and we write $\norm{f}_{B_{\ell^n_{p}}}=\sup_{z\in B_{\ell^n_{p}}}\norm{f(z)}_{X}$.\\[2mm]
\indent A bounded linear operator $v:X\rightarrow Y$ between two Banach spaces $X$ and $Y$ is called $(p,q)$-summing, $1\leq p,q<\infty$, if there is a constant $C>0$ such that for each choice of finitely many vectors $x_1,\ldots,x_n\in X$, we have
\bea\label{b5}
	\left(\sum_{k=1}^{n}\norm{v(x_k)}^p\right)^{\frac{1}{p}}\leq C\sup_{x^*\in B_{X^*}}\left(\sum_{k=1}^{n}|x^*(x_k)|^q\right)^{\frac{1}{q}}.
\eea
The best such constant $C$ satisfying \eqref{b5} is denoted by $\pi_{p,q}(v)$. In particular, if $p=q$, we say that $v$ is $p$-summing and we write the best constant as $\pi_p(v)$.\\[2mm]
\indent In this paper, we mainly focus on the unit ball 
\beas B_{\ell^n _p}:= \left\{z =(z_1, z_2,\ldots, z_n)\in \Bbb{C}^n:\norm{z}_{p}<1\right\}\eeas
of Minkowski space $\ell^n_p$ $(1\leq p\leq \infty)$, where the $p$-norm is given by $\norm{z}_{p}=\left(\sum_{k=1}^{n}|z_k|^p\right)^{1/p}$ for $1\leq p<\infty$ and $\norm{z}_\infty=\sup_{1\leq k\leq n} |z_k|$ for $p=\infty$.\\[2mm]
\indent  
Before proceeding further, let us first introduce the following definition of multidimensional Bohr radius for holomorphic functions defined on $B_{\ell^n_p}$ $( 1\leq p\leq \infty)$ with values in a complex Banach space $X$, primarily motivated by the work of Defant {\it et al.} \cite{DMS2012} and Galicer {\it et al.} \cite{GMM2020}.
\begin{defi}\label{Defi1}
Let $1\leq p,q\leq \infty$, $1\leq s<\infty$, $n\in \Bbb{N}$, and let $v: X \rightarrow Y$ be a bounded linear operator between complex Banach spaces $X$ 
and $Y$ such that $\norm{v}\leq \lambda^{1/s}$. The $\lambda_{s}$-mixed Bohr radius of $v$, denoted by $K^s(B_{\ell^n_{p}},B_{\ell^n_{q}}, v,\lambda)$, is the supremum of all $r\geq 0$ such that the following inequality 
	\bea\label{a1}
		\sup_{z\in rB_{\ell^n_{q}}}\,\,\sum_{\alpha\in \Bbb{N}^n_{0}}\norm{v(c_{\alpha})z^{\alpha}}^{s}_{Y}\leq \lambda\sup_{z\in B_{\ell^n_{p}}}\norm{\sum_{\alpha\in \Bbb{N}^n_{0}}c_{\alpha}z^{\alpha}}^s_{X}
	\eea
holds for all holomorphic functions $f(z)=\sum_{\alpha\in \Bbb{N}^n_{0}}c_{\alpha}z^{\alpha}$ on $B_{\ell^n_{p}}$ with values in $X$. 
We write $K^s(B_{\ell^n_p}, B_{\ell^n_q},v)$ whenever $\lambda=1$. If $v$ is an identity operator on $X$, then we simply write $K^s(B_{\ell^n_p},B_{\ell^n_q},X,\lambda)$ and $K^s(B_{\ell^n_p},X,\lambda)$ whenever $p=q$. For $X=\Bbb{C}$, we write $K^s(B_{\ell^n_p},B_{\ell^n_q},\lambda)$ for $K^s(B_{\ell^n_p},B_{\ell^n_q},\Bbb{C},\lambda)$ and $K(B_{\ell^n_p},\lambda)$ whenever $p=q$ and $s=1$.
\end{defi}
Evidently, the inequality (\ref{a1}) holds trivially if $f$ is an unbounded holomorphic function on $B_{\ell^n_{p}}$. Thus, throughout the paper, we assume that the holomorphic 
function $f$ is bounded on $B_{\ell^n_{p}}$. Clearly, using the above notation, \textrm{Theorem A} reduces to $K(B_{\ell^1_p})=1/3$, where $\Bbb{D}=B_{\ell^1_p}$. 
Again, if $v$ is a null operator, then this supremum is $\infty$. Thus, to avoid the trivial case, we assume throughout the paper that $v\neq 0$. Furthermore, we observe that 
\bea\label{a3}
K^s(B_{\ell^n_p},B_{\ell^n_q},v,\lambda)\geq \max \left\{K^s\left(B_{\ell^n_p},B_{\ell^n_q},X,\lambda/\norm{v}^s\right), K^s\left(B_{\ell^n_p},B_{\ell^n_q},Y,\lambda/\norm{v}^s\right)\right\}. \eea
Note that for two sequences of positive numbers $(a_n)$ and $(b_n)$, we write $a_n\prec b_n$ if there exists a
constant $\delta > 0$ such that $a_n\leq \delta b_n$ for every $n$, and $a_n\sim b_n$ whenever $a_n\prec b_n$ and $b_n\prec a_n$.
 In the classical sense, Bombieri and Bourgain \cite{B1962,BB2004} have studied the general Bohr radius $K(\Bbb{D},\lambda)$ whenever $\lambda\geq 1$. In fact, Bombieri and Bourgain \cite{B1962,BB2004} have derived the following exact value 
\beas
	K(B_{\ell^1_p},\lambda)=\frac{1}{3\lambda-2\sqrt{2(\lambda^2-1)}}\quad\text{for all}\quad 1\leq \lambda\leq \sqrt{2},
\eeas
and $K(B_{\ell^1_p},\lambda)\sim \sqrt{\lambda^2-1}/\lambda$ for $\lambda \rightarrow \infty$. Let $H^\infty(\mathbb{D}, X)$ denote the space of all bounded analytic 
functions from the unit disk $\mathbb{D}$ into a complex Banach space $X$, equipped with 
the norm $\norm{f}_{H^\infty(\mathbb{D}, X)} := \sup_{|z|<1} \norm{f(z)}_X$.
In 2010, Blasco \cite{B2010} defined the Bohr radius $R(X)$ of $X$ by 
\be \label{c3}R(X)=\sup\left\{r\geq 0: \sum_{k=0}^{\infty} \norm{ x_k}  r^k\leq \norm{ f}_{H^\infty(\mathbb{D}, X)}\ \text{for every }\ f \in H^\infty(\mathbb{D}, X),\ f(z) = \sum_{k=0}^{\infty} x_k z^k\right\},\ee
and established that for $X=\ell^2_p$, we have $R(X)=0$, {\it i.e.}, $K^1(\Bbb{D},\ell^2_p,1)=0$ for $1\leq p\leq \infty$ (in the notation of 
Definition \ref{Defi1}).
Thus, the study of Bohr radius problem becomes irrelevant for $\text{dim}(X)>1$. 
Hence, it seems reasonable to modify the Bohr inequality to obtain more interesting results. This fact might have motivated Defant {\it et al.} \cite{DMS2012} to modify the definition of the
Bohr radius (\ref{c3}) in order to obtain non-trivial Bohr radii.\\[2mm]
\indent The study of multi-dimensional Bohr radii regained attention through the research conducted by Dineen and Timoney \cite{DT1989}. Motivated by this work, Boas and Khavinson \cite{BK1997} showed that 
\beas \frac{1}{3\sqrt{n}}\leq K_n\leq 2\sqrt{\frac{\log n}{n}}\quad \text{for}\quad n\geq 2.\eeas
Later, Boas \cite{B2000} has established the bound of $K(B_{\ell^n_p})$ for arbitrary $p$ and also conjectured that {\it ``This logarithmic factor, an artifact of the proof, presumably 
should not really be present''}. For some time it was thought that the lower bound of $K_n$ could not be improved. In 2006, Defant and Frerick \cite{DF2006} disproved this 
conjecture by obtaining a logarithmic lower bound which is almost correct asymptotic estimate for the Bohr radius $K(B_{\ell^n_p})$ for $1\leq p\leq \infty$. 
Using the link between the multi-dimensional Bohr radii and local Banach space theory (see \cite{DGM2003}), Defant and Frerick \cite{DF2006} have proved that there is a 
constant $C>0$ such that for every $1\leq p\leq \infty$ and $n$, we have
 \beas
 	\frac{1}{C}\left(\frac{\log n}{n\log(\log n)}\right)^{1-\frac{1}{\min\{p,2\}}}\leq K(B_{\ell^n_p}).
 \eeas
 In view of results of Aizenberg \cite{A2000}, Boas and Khavinson \cite{BK1997}, Boas \cite{B2000}, Dineen and Timoney \cite{DT1989}, and Defant and Frerick \cite{DF2006}, there is a constant $C\geq1$ such that for all $1\leq p\leq \infty$ and all $n\in\Bbb{N}$, we have
 \beas
 	\frac{1}{C}\left(\frac{\log n}{n\log(\log n)}\right)^{1-\frac{1}{\min\{p,2\}}}\leq K(B_{\ell^n_p})\leq C\left(\frac{\log n}{n}\right)^{1-\frac{1}{\min\{p,2\}}}.
 \eeas
 Later, Defant and Frerick \cite{DF2011} made a significant improvement on the above result, showing that there exists a constant $C \geq 1$ such that for all $1\leq p\leq \infty$ and all $n\in\Bbb{N}$, we have
 \beas
 	\frac{1}{C}\left(\frac{\log n}{n}\right)^{1-\frac{1}{\min\{p,2\}}}\leq K(B_{\ell^n_p}).
 \eeas
The case $p=\infty$ was settled by Defant {\it et al.} \cite{DFOOS2011} and they showed that $K(B_{\ell^n_{\infty}})$ behaves asymptotically as $\sqrt{\log n/n}$ as a consequence of the fact that Bohnenblust-Hille inequality for homogeneous polynomials is hypercontractive. Furthermore, Bayart {\it et al.} \cite{BPS2014} proved that
 \beas
 	\lim_{n\rightarrow \infty}\frac{K(B_{\ell^n_{\infty}})}{\sqrt{\frac{\log n}{n}}}=1.
 \eeas
 The gap between the upper and lower bounds in all the aforementioned papers motivated many researchers to find the exact asymptotic order of $K(B_{\ell^n_p})$ for $1\leq p\leq \infty$. The Bohr radius has been studied extensively in several directions such as complex analysis \cite{BDS2019, DMP2009}, number theory \cite{BB2004, DFOOS2011} and operator algebra \cite{PPS2002}. In 2012, Defant {\it et al.}\cite{DMS2012} estimated the asymptotic behavior of $K(B_{\ell^n_{\infty}}, X,\lambda)$ for the finite and infinite dimensional Banach space $X$. Later, Galicer {\it et al.} \cite{GMM2020} determined the exact asymptotic growth of the mixed $(p,q)$-Bohr radius as $n\rightarrow\infty$ for every $1\leq p,q\leq\infty$. \\[2mm]
\noindent 
 As demonstrated by the scalar case, studying the $m$-homogeneous polynomial first is typically a beneficial approach for obtaining nontrivial estimates for Bohr radii. The following definition is the $m$-homogeneous analogue of Definition \ref{Defi1} for $m\in \Bbb{N}$.
\begin{defi}\label{Defi2}
	Let $1\leq p,q\leq \infty$, $1\leq s<\infty$, $m,n \in \Bbb{N}$, and $v: X\rightarrow Y$ be a bounded operator in complex Banach spaces such that $\norm{v}\leq \lambda^{1/s}$. Then $K^s_m\left(B_{\ell^n_p},B_{\ell^n_q},v,\lambda\right)$ is defined to be the supremum of all $r\geq 0$ such that the following inequality 
	\beas
		\sup_{z\in rB_{\ell^n_{q}}}\ \sum_{\alpha \in\Lambda(m,n)}\norm{v(c_{\alpha})z^{\alpha}}^{s}_{Y}\leq \lambda\sup_{z\in B_{\ell^n_{p}}}\norm{\sum_{\alpha \in\Lambda(m,n)}c_{\alpha}z^{\alpha}}^s_{X}.
	\eeas
holds	for every $m$-homogeneous polynomials $P\in \mathcal{P}(^m\ell^n_p, X)$ with $P(z)=\sum_{\alpha \in\Lambda(m,n)}c_{\alpha}z^{\alpha}$. 
\end{defi}
The constant $K^s_m\left(B_{\ell^n_p},B_{\ell^n_q},v,\lambda\right)$ in Definition \ref{Defi2} can be rewritten as
\begin{small}\be\label{Defi3}
	K^s_m\left(B_{\ell^n_p},B_{\ell^n_q},v,\lambda\right)= \sup \left\{r> 0 : \sup_{z\in B_{\ell^n_{q}}}\ \sum_{\alpha \in\Lambda(m,n)}\norm{v(c_{\alpha})z^{\alpha}}^s_Y\leq \frac{\lambda}{r^{ms}}\norm{P}_{\mathcal{P}(^m\ell^n_p,X)}^s\,\,\text{for} \,\, P\in \mathcal{P}(^m\ell^n_p, X)\right\}.
\ee\end{small}
Moreover, it is easy to observe that 
\bea\label{a2}
	K^s_m\left(B_{\ell^n_p},B_{\ell^n_q},v,\lambda\right)=\sqrt[ms]{\lambda}K^s_m\left(B_{\ell^n_p},B_{\ell^n_q},v\right).
\eea
The powered versions of the Bohr radii $K^p(\Bbb{D})$ and $K^p(B_{\ell^n_{\infty}})$ for scalar-valued holomorphic functions were studied by Djakov and Ramanujan 
\cite{DR2000}. In 2004, B\'{e}n\'{e}teau {\it et al.} \cite{BDK2004} explored the analogous radii in the context of Hardy spaces. Apart from the scalar valued functions, Blasco 
\cite{B2017} has investigated the radius $K^p(\Bbb{D},X)$ for arbitrary complex Banach spaces. Recently, Das \cite{D2023,D2024} has obtained the exact asymptotic value of 
$K^p(B_{\ell^n_{\infty}})$ and investigated $p$-Bohr radius for Hardy space of $X$-valued holomorphic functions defined on complete Reinhardt domain. Recently, Allu and Pal 
\cite{AP2026} studied multidimensional Bohr radii for Banach space valued holomorphic functions and obtained sharp asymptotic estimates for bounded linear operators between complex Banach spaces in both finite and infinite dimensions. Their work provides a natural motivation for the present paper.\\[2mm]
 \indent To explore more general aspects of multidimensional Bohr radii, we introduce $\lambda,p,q$ and $s$ in our Definition \ref{Defi1}, in more general sense. The aim of this 
 paper is to obtain estimates for the Bohr radii $K^s(B_{\ell^n_p},B_{\ell^n_q},v,\lambda)$ for various choices of the bounded linear operator $v$ between complex Banach spaces. We first prove 
 that $K^s(B_{\ell^n_p},B_{\ell^n_q},v,\lambda)>0$ for every bounded operator $v$ with $\norm{v}<\lambda^{1/s}$. For an infinite-dimensional Banach space $X$, we obtain 
 the exact asymptotic behavior of $K^s(B_{\ell^n_p},B_{\ell^n_q},v,\lambda)$. Finally, we study $K^s(B_{\ell^n_p},B_{\ell^n_q},v,\lambda)$ 
 for a Banach space of cotype $r$ and for an $r$-concave Banach lattice.\\[2mm]
\noindent \textbf{Novelty and significance.} In this paper, we introduce and study a generalized framework for multidimensional mixed Bohr radii of Banach space valued holomorphic 
functions defined on the unit ball of $\ell_p^n$, where the coefficients are acted upon by a bounded linear operator $v: X \to Y$. The main novelty of this paper is the derivation of 
positivity and dimension-dependent estimates for these mixed 
Bohr radii, relating their behavior to geometric properties of the spaces, including cotype and concavity, as well as to summing properties of the operators. In particular, 
we obtain a lower bound involving the factor $n^{-\max\{0,\,1/p-1/q\}}((\log n)/n)^{1-1/r}$ for concave Banach lattices and summing 
operators, and recovers the classical Defant--Frerick exponent in the scalar case. These results connect classical scalar Bohr phenomena with operator valued and Banach 
space valued settings, and highlight how the geometry of the spaces and the mapping properties of the operator influence coefficient estimates and their dependence on dimension.
\section{Key Lemmas}
The following lemma provides a connection between Definitions \ref{Defi1} and \ref{Defi2}, which will be one of the most important tools to prove our main results in this paper.
\begin{lem}\label{lem1}
	Let  $1\leq p,q\leq \infty$, $1\leq s<\infty$ and $v:X\rightarrow Y$ be a bounded operator between complex Banach spaces $X$ and $Y$ with $\norm{v}<\lambda^{1/s}$. Then, we have 
	\begin{enumerate}
		\item[(i)] $\left(\dfrac{\lambda-\norm{v}^s}{2\lambda-\norm{v}^s}\right)^{1/s}\inf\limits_{m\in \Bbb{N}}K^s_m(B_{\ell^n_p},B_{\ell^n_q}, v, \lambda)\leq 
		K^s(B_{\ell^n_p},B_{\ell^n_q}, v, \lambda)\leq \inf\limits_{m\in \Bbb{N}}K^s_m(B_{\ell^n_p},B_{\ell^n_q}, v, \lambda)$,\\[2mm]
		\item[(ii)] $\left(\dfrac{\lambda-\norm{v}^s}{\lambda-\norm{v}^s+1}\right)^{1/s}\inf\limits_{m\in\Bbb{N}}K^s_m(B_{\ell^n_p},B_{\ell^n_q}, v)\leq 
		K^s(B_{\ell^n_p},B_{\ell^n_q}, v, \lambda)$\\[4mm] and $K^s(B_{\ell^n_p},B_{\ell^n_q}, v, \lambda)\leq \max\left\{1,\lambda^{1/s}\right\} \inf\limits_{m\in \Bbb{N}}K^s_m(B_{\ell^n_p},B_{\ell^n_q}, v)$.
	\end{enumerate}
\end{lem}
\begin{proof}
(i) Since $K^s(B_{\ell^n_p},B_{\ell^n_q}, v, \lambda)\leq K^s_m(B_{\ell^n_p},B_{\ell^n_q}, v, \lambda)$ for all $m\geq 1$, the right hand side inequality of (i) is obvious.  
Let $f(z)=\sum_{\alpha\in \Bbb{N}^n_{0}}c_{\alpha}z^{\alpha}$ be a holomorphic function on $B_{\ell^n_p}$ with values in $X$ and $x\in \Bbb{C}^n$ such that 
$\norm{x}_q\leq \inf\limits_{m\in \Bbb{N}} K^s_m(B_{\ell^n_p},B_{\ell^n_q}, v, \lambda)$.  Since $\norm{x}_q\leq  K^s_m(B_{\ell^n_p},B_{\ell^n_q}, v, \lambda)$ for all $m\in \Bbb{N}$, we have 
\bea\label{a4}  \sum_{\alpha \in\Lambda(m,n)}\norm{v(c_{\alpha}) x^\alpha}^s\leq \lambda\sup_{z\in B_{\ell^n_{p}}}\norm{P_m(z)}^s_{X},\eea
where  $P_m(z)=\sum_{\alpha \in\Lambda(m,n)}c_{\alpha}z^{\alpha}$ is the $m$-homogeneous part of $f$. 
For each fixed $z\in B_{\ell_p^n}$, we define the function $g_z:\{\lambda\in\mathbb C:|\lambda|<1/\|z\|_p\}\to X$ by  $g_z(\lambda)=f(\lambda z)$ when $z\neq 0$, and 
$g_0(\lambda)=f(0)$. 
For $z \neq 0$, the disk $\{\lambda : |\lambda| < 1/\norm{z}_p\}$ contains the closed unit disk $\overline{\mathbb{D}}\) in its interior and 
$g_z(\lambda) = \sum_{m=0}^{\infty} \lambda^m P_m(z)$ is holomorphic on this disk. For $z = 0$, the function $g_0(\lambda) = f(0)$ is constant and hence holomorphic on 
$\Bbb{C}$ with $g_0(\lambda)= \sum_{m=0}^{\infty} \lambda^m P_m(0)$, where $P_m(0) = 0$ for $m \geq 1$. In both cases, $g_z$ is holomorphic on a
disk containing $\overline{\mathbb{D}}$ in its interior.
By the one-variable Cauchy coefficient formula, we have
\beas
P_m(z)=\frac{1}{2\pi}\int_0^{2\pi} f(e^{i\theta}z)e^{-im\theta}\,d\theta\quad \mbox{for every $m \in\Bbb{N}\cup\{0\}$ and every $z\in B_{\ell_p^n}$}.
\eeas
Since $B_{\ell_p^n}$ is balanced, it follows that $e^{i\theta}z\in B_{\ell_p^n}$ for every $\theta\in\Bbb{R}$ and consequently, $\norm{f(e^{i\theta}z)}_X\leq \norm{f}_{B_{\ell_p^n}}$. Therefore,
\beas \norm{P_m(z)}_X\leq \frac{1}{2\pi}\int_0^{2\pi}\norm{f(e^{i\theta}z)}_X\,d\theta\leq \norm{f}_{B_{\ell_p^n}}.\eeas
Taking the supremum over $z\in B_{\ell_p^n}$ gives $\norm{P_m}_{\mathcal P({}^m\ell_p^n,X)} \leq \norm{f}_{B_{\ell_p^n}}$.
It is evident that for any $t \in(0, 1)$, we have $t\,x\in K^s_m(B_{\ell^n_p},B_{\ell^n_q}, v, \lambda) B_{\ell^n_{q}}$ for all $m\in \Bbb{N}$. 
	Therefore, we obtain
	\beas
		\sum_{\alpha\in \Bbb{N}^n_{0}}\norm{v(c_{\alpha})\left(t \,x\right)^{\alpha}}^s_{Y} &=& \norm{v(c_0)}^s+\sum_{m=1}^{\infty}t^{ms} \left(\sum_{\alpha \in\Lambda(m,n)}\norm{v(c_{\alpha})x^\alpha}^s\right) \\[2mm]
		&\leq& \norm{v}^s\norm{c_0}^s + \sum_{m=1}^{\infty} \lambda t^{ms}  \norm{P_m}^s_{\mathcal{P}(^m\ell^n_p,X)}\\[2mm]
		&\leq& \left(\norm{v}^s + \lambda \sum_{m=1}^{\infty} t^{ms}\right)\norm{f}^s_{B_{\ell_p^n}}\\[2mm]
		&=& \left(\norm{v}^s + \lambda \cdot\frac{t^s}{1-t^s}\right)\norm{f}^s_{B_{\ell_p^n}}.
	\eeas
	If we choose $t=\left(\lambda-\norm{v}^s\right)^{1/s}/\left(2\lambda-\norm{v}^s\right)^{1/s}$, then we observe that $t\in (0,1)$ and  
	\beas
		\sum_{\alpha\in \Bbb{N}^n_{0}}\norm{v(c_{\alpha})\left(t \,x\right)^{\alpha}}^s_{Y} \leq \lambda \norm{f}^s_{B_{\ell_p^n}}.
	\eeas
	Thus, for any $x\in \Bbb{C}^n$ with $\norm{x}_q\leq \inf\limits_{m\in \Bbb{N}}\{K^s_m(B_{\ell^n_p},B_{\ell^n_q}, v, \lambda)\}$, the point $t\,x$ satisfies the $\lambda_s$-mixed Bohr inequality of $v$. Hence, we have 
	\beas \left(\frac{\lambda-\norm{v}^s}{2\lambda-\norm{v}^s}\right)^{1/s}\inf\limits_{m\in \Bbb{N}}K^s_m(B_{\ell^n_p},B_{\ell^n_q}, v, \lambda)\leq 
		K^s(B_{\ell^n_p},B_{\ell^n_q}, v, \lambda).\eeas
	(ii) From (\ref{a2}), we have $K^s_m\left(B_{\ell^n_p},B_{\ell^n_q},v,\lambda\right)=\lambda^{1/(ms)}K^s_m\left(B_{\ell^n_p},B_{\ell^n_q},v\right)$. Using the fact that 
	$K^s(B_{\ell^n_p},B_{\ell^n_q}, v,\lambda)\leq \inf\limits_{m\in\Bbb{N}} K^s_m(B_{\ell^n_p},B_{\ell^n_q}, v,\lambda)$, we obtain
	\beas
		K^s\left(B_{\ell^n_p}, B_{\ell^n_q},v,\lambda\right)\leq \inf\limits_{m\in \Bbb{N}} \left\{\lambda^{1/(ms)}K^s_m\left(B_{\ell^n_p}, B_{\ell^n_q},v\right)\right\}\leq \max\left\{1,\lambda^{1/s}\right\} \inf\limits_{m\in \Bbb{N}} K^s_m\left(B_{\ell^n_p}, B_{\ell^n_q},v\right).
	\eeas
	For the left hand side inequality, we suppose again that $x\in \Bbb{C}^n$ with $\norm{x}_q\leq \inf\limits_{m\in \Bbb{N}}  K^s_m(B_{\ell^n_p}, B_{\ell^n_q},v)$. Hence, we have 
	\beas\sum_{\alpha \in\Lambda(m,n)}\norm{v(c_{\alpha}) x^\alpha}^s\leq \sup_{z\in B_{\ell^n_{p}}}\norm{P_m(z)}^s_{X},\eeas
where  $P_m(z)=\sum_{\alpha \in\Lambda(m,n)}c_{\alpha}z^{\alpha}$.
	Thus, for any $t\in(0, 1)$, we have $t\,x\in K^s_m(B_{\ell^n_p}, B_{\ell^n_q},v)B_{\ell^n_{q}}$, and consequently, we obtain
	\beas
		\sum_{\alpha\in \Bbb{N}^n_{0}} \norm{v(c_{\alpha})\left(t\, x\right)^{\alpha}}^s_{Y}
		&=&\norm{v(c_0)}^s+\sum_{m=1}^{\infty}t^{ms} \sum_{\alpha \in\Lambda(m,n)}\norm{v(c_{\alpha})x^\alpha}^s\\[2mm]
		&\leq& \norm{v}^s\norm{c_0}^s + \sum_{m=1}^{\infty}t^{ms} \norm{P_m}^s_{\mathcal{P}(^m\ell^n_p,X)}\\[2mm]
		&\leq& \left(\norm{v}^s +  \sum_{m=1}^{\infty} t^{ms}\right)\norm{f}^s_{B_{\ell^n_p}}\\[2mm]
		&=& \left(\norm{v}^s + \frac{t^s}{1-t^s}\right)\norm{f}^s_{B_{\ell^n_p}}.
	\eeas
	If we choose $t=\left(\lambda-\norm{v}^s\right)^{1/s}/\left(\lambda-\norm{v}^s+1\right)^{1/s}$, then we observe that $t\in (0,1)$ and  
	\beas
		\sum_{\alpha\in \Bbb{N}^n_{0}}\norm{v(c_{\alpha})\left(t \,x\right)^{\alpha}}^s_{Y} \leq \lambda\norm{f}^s_{B_{\ell^n_p}}.
	\eeas
	Hence, we have 
	\beas \left(\frac{\lambda-\norm{v}^s}{\lambda-\norm{v}^s+1}\right)^{1/s}\inf\limits_{m\in\Bbb{N}}K^s_m(B_{\ell^n_p},B_{\ell^n_q}, v)\leq 
		K^s(B_{\ell^n_p},B_{\ell^n_q}, v, \lambda).\eeas
	This completes the proof.
\end{proof}

\begin{rem}
Let $v:X\rightarrow Y$ be a bounded operator between complex Banach spaces $X$ and $Y$ with $\norm{v}<\lambda^{1/s}$. Define $ r_0 := \inf\limits_{m \in \Bbb{N}} K^s_m(B_{\ell^n_p},B_{\ell^n_q}, v, \lambda)$.
Then for every bounded holomorphic function $f: B_{\ell^n_p} \rightarrow X$ with $f(0) = 0$, we have
\beas 
		\sup_{z\in 2^{-1/s} r_0  B_{\ell^n_{q}}}\ \sum_{\alpha(\not=0) \in \Bbb{N}^n_0}\norm{v(c_{\alpha})z^{\alpha}}^{s}_{Y}\leq \lambda\sup_{z\in B_{\ell^n_{p}}}\norm{\sum_{\alpha(\not=0) \in \Bbb{N}^n_0}c_{\alpha}z^{\alpha}}^s_{X}.
\eeas
Consequently,
\beas \frac{1}{2^{1/s}}\inf\limits_{m\in \Bbb{N}}K^s_m(B_{\ell^n_p},B_{\ell^n_q}, v, \lambda)\leq 
		K^s(B_{\ell^n_p},B_{\ell^n_q}, v, \lambda).\eeas
For bounded holomorphic functions vanishing at the origin, the factor $\left(\lambda-\norm{v}^s\right)^{1/s}/\left(2\lambda-\norm{v}^s\right)^{1/s}$ in Lemma~\ref{lem1} 
can be replaced by the universal constant $2^{-1/s}$, independently of $\lambda$ and $\norm{v}$.
\end{rem}
\begin{lem}\cite[Theorem 3.2]{BPV2004}\label{lem2}
Let $X_1, \ldots, X_m$  and $Y $ be complex Banach spaces, and let $T: X_1 \times \cdots \times X_m \to Y $ be a bounded $m$-linear operator. 
Suppose $Y$ has cotype $r$ with $2 \leq r < \infty$. Then, $T$ is multiple $(r,1)$-summin and
\beas \pi_{(r,1)}(T) \leq C_r(Y)^m \norm{T},\eeas
where $C_r(Y)$ is the cotype $r$ constant of $Y$. Equivalently, for every 
choice of finitely many vectors $x_{i_j}^j \in X_j$ $(1 \leq i_j \leq n_j , 1 \leq j \leq m )$, we have
\beas
\left(\sum_{i_1, \ldots, i_m = 1}^{n_1, \ldots, n_m} \norm{T(x_{i_1}^1, \ldots, x_{i_m}^m)}^r\right)^{1/r}\leq C_r(Y)^m \norm{T}\prod_{j=1}^m\norm{(x_{i_j}^j)_{i_j=1}^{n_j} }_{\ell_1^w(X_j)},\eeas
where $\norm{(x_{i_j}^j)_{i_j=1}^{n_j}}_{\ell_1^w(X_j)}$ denotes the weak $\ell_1$-norm of the sequence $(x_{i_j}^j)_{i_j=1}^{n_j}$ in $X_j$, defined by
\beas
\left\| (x_{i_j}^j)_{i_j=1}^{n_j} \right\|_{\ell_1^w(X_j)}= \sup_{x^* \in B_{X_j^*}} \sum_{i_j=1}^{n_j} \left|x^*\left(x_{i_j}^j \right)\right|.\eeas
\end{lem}
\noindent The following lemma provides an explicit bound for the coefficients of a symmetric $m$-linear operator on $\ell_p^n$ in terms of its operator norm and the dimension $n$.
\begin{lem}\label{lem3}
Let $X$ be a complex Banach space of cotype $r$ with $2 \leq r < \infty$ and let $A: \underbrace{\ell_p^n\times \cdots\times \ell_p^n}_{m\text{-times}} \to X$ be a bounded 
symmetric $m$-linear operator with operator norm
\beas
\|A\| = \sup_{z^{(j)} \in B_{\ell_p^n}} \norm{A(z^{(1)}, \ldots, z^{(m)})}_X.\eeas
Let $(e_{i_j})_{i_j=1}^n=(e_1,\ldots,e_n)$ be the canonical basis of $\ell_p^n$, where $j=1,\ldots,m$. Then
\beas \left(\sum_{i\in \mathcal{M}(m,n)}\norm{A(e_{i_1}, \ldots, e_{i_m})}^r\right)^{1/r}\leq C_r(X)^m \, n^{m/p} \, \norm{A},\eeas
where $C_r(X)$ is the cotype $r$ constant of $X$ and $n^{1/\infty} = 1$.
\end{lem}
\begin{proof}
Let $X_1 = \cdots = X_m = \ell_p^n$, $Y = X$,  $T = A$. For each $j = 1,\dots,m$ and $i_j = 1,\dots,n$, let $x_{i_j}^j = e_{i_j}$, where $e_{i_j}$ is the $i_j$-th canonical basis vector of $\ell_p^n$, assigned to the $j$-th variable. In view of \textrm{Lemma \ref{lem2}}, we obtain
\bea\label{c1}
\left(\sum_{i\in \mathcal{M}(m,n)}\norm{A(e_{i_1}, \ldots, e_{i_m})}^r\right)^{1/r}\leq C_r(X)^m \norm{A}\prod_{j=1}^m\norm{ (e_{i_j})_{i_j=1}^n}_{\ell_1^w(\ell_p^n)}.
\eea
For the canonical basis $(e_j)_{j=1}^n$ of $\ell_p^n$, the weak $\ell_1$-norm 
is \beas \|(e_j)_{j=1}^n\|_{\omega_1(\ell_p^n)} = \sup_{x^* \in B_{\left(\ell_{p}^n\right)^*}} \sum_{j=1}^n \left|x^*\left(e_j \right)\right| = \sup_{x^* \in B_{\left(\ell_{p}^n\right)^*}} \sum_{j=1}^n \left|x_j^*\right|,\eeas 
where $(x_j^*)_{j=1}^n\in B_{\ell_{p'}^n}$ with $1/p + 1/p' = 1$. By H\"{o}lder's inequality, we have 
\beas \sum_{j=1}^n |x_j^*| \leq  \left(\sum_{j=1}^n |x_j^*|^{p'}\right)^{1/{p'}}\left(\sum_{j=1}^n 1\right)^{1/p} \leq  n^{1/p}.\eeas
Therefore, we have $\|(e_j)_{j=1}^n\|_{\omega_1(\ell_p^n)}=n^{1/p}$ and 
\bea\label{c2} \prod_{j=1}^m \|(e_{i_j})_{i_j=1}^n\|_{\omega_1(\ell_p^n)} = (n^{1/p})^m = n^{m/p}.\eea
From (\ref{c1}) and (\ref{c2}), we obtain the desired result. This completes the proof.
\end{proof}
\section{Main Results}\label{Pal-Vasu-P5-section-03}
\noindent In the following result, we prove that the $\lambda_{s}$-mixed Bohr radius $K^s(B_{\ell^n_p}, B_{\ell^n_q},v,\lambda)$ of $v$ is positive.
\begin{thm}\label{T1}
 	Let $1\leq p, q\leq \infty$, $1\leq s<\infty$ and $v:X\rightarrow Y$ be a bounded non-null operator between complex Banach spaces $X$ and $Y$ such that $\norm{v}<\lambda^{1/s}$. Then there exists a constant $C>0$ such that for all $n\in \Bbb{N}$, we have 
 	\beas
 		K^s(B_{\ell^n_p}, B_{\ell^n_q},v,\lambda)\geq C \frac{1}{e\, n^{\max\left\{1, 1+\frac{1}{s}-\frac{1}{q}\right\}}},
 	\eeas
 	where 
 	\beas 
 	C=\left\{\begin{array}{ll}
 		\max \left\{\left(\dfrac{\lambda-\norm{v}^{s}}{2\lambda-\norm{v}^s}\right)^{1/s}, \dfrac{1}{\norm{v}}\left(\dfrac{\lambda-\norm{v}^s}{\lambda-\norm{v}^s+1}\right)^{1/s}\right\}, & \mbox{ if  \,$\norm{v}\geq 1$},\\[5mm] 
 		\max \left\{\left(\dfrac{\lambda-\norm{v}^s}{2\lambda-\norm{v}^s}\right)^{1/s}, \left(\dfrac{\lambda-\norm{v}^s}{\lambda-\norm{v}^s+1}\right)^{1/s}\right\}, & \mbox{ if \,$0<\norm{v}<1$}.
 	\end{array}\right.
 	\eeas
 \end{thm}
 \begin{proof}
	Let $P(z)=\sum_{\alpha \in\Lambda(m,n)}c_{\alpha}z^{\alpha}\in \mathcal{P}(^m \ell^n_p,X)$ be an $m$-homogeneous polynomial. For any $\varphi\in X^*$ with $\norm{\varphi}\leq 1$, the scalar function $g=\varphi \circ P: B_{\ell^n_p}\rightarrow\Bbb{C}$ is an $m$-homogeneous polynomial with expansion 
	\beas g(z)=\sum_{\alpha \in\Lambda(m,n)}\varphi\left(c_{\alpha}\right)z^{\alpha}.\eeas
	Applying the Cauchy's estimate to $g$ gives
	\bea\label{c4}\left|\varphi(c_\alpha)\right|\leq \left(\frac{|\alpha|^{|\alpha|}}{\alpha^\alpha} \right)^{1/p}\norm{g}_{\mathcal{P}(^m \ell_p^n, \Bbb{C})}\leq \frac{|\alpha|^{|\alpha|}}{\alpha^\alpha}\norm{g}_{\mathcal{P}(^m \ell_p^n, \Bbb{C})},\eea
	where $\norm{g}_{\mathcal{P}(^m \ell_p^n, \Bbb{C})}=\sup_{z\in B_{\ell_p^n}}|\varphi(P(z))|$ (see \cite{B2000, DF2006}). Since $|\varphi(P(z))|\leq \norm{\varphi} \norm{P(z)}_X$, we have 
	\beas \left|\varphi(c_\alpha)\right|\leq \frac{|\alpha|^{|\alpha|}}{\alpha^\alpha}\norm{\varphi} \sup_{z\in B_{\ell_p^n}}\norm{P(z)}_X.\eeas
	Taking supremum over all $\varphi\in B_{X^*}$ and applying the Hahn-Banach theorem in (\ref{c4}), we obtain
		\bea\label{a5} \norm{c_\alpha}_X=\sup_{\varphi\in B_{X^*}}\left|\varphi(c_\alpha)\right|\leq \frac{|\alpha|^{|\alpha|}}{\alpha^\alpha}\norm{P}_{\mathcal{P}(^m \ell_p^n, X)}.\eea
	Using (\ref{a5}), we obtain
	\bea\label{a6}
		\sup_{z\in B_{\ell_q^n}}\sum_{\alpha \in\Lambda(m,n)}\norm{v(c_{\alpha})z^{\alpha}}^s &\leq&  \sup_{z\in B_{\ell_q^n}}\sum_{\alpha \in\Lambda(m,n)}\norm{v}^s \norm{c_\alpha}_X^s|z^{\alpha}|^s \nonumber\\[2mm]
		& \leq &\norm{v}^s \norm{P}_{\mathcal{P}(^m \ell_p^n, X)}^s  \sup_{z\in B_{\ell_q^n}}\sum_{\alpha \in\Lambda(m,n)}\left(\frac{m^m}{\alpha^\alpha}\right)^s |z^{\alpha}|^s. \eea
	We observe that 
	\bea\label{a7} \frac{m^m}{m!}\leq e^m\leq \frac{\alpha^\alpha}{\alpha!}e^m\quad \text{and}\quad 1\leq \frac{m!}{\alpha!}\leq n^m.\eea
	Using (\ref{a7}) in (\ref{a6}), we obtain 	
		\beas \sup_{z\in B_{\ell_q^n}}\sum_{\alpha \in\Lambda(m,n)}\norm{v(c_{\alpha})z^{\alpha}}^s & \leq &\norm{v}^s \norm{P}_{\mathcal{P}(^m \ell_p^n, X)}^s (en)^{ms} \sup_{z\in B_{\ell_q^n}}\sum_{\alpha \in\Lambda(m,n)} |z^{\alpha}|^s\\[2mm]
		&\leq &\norm{v}^s \norm{P}_{\mathcal{P}(^m \ell_p^n, X)}^s (en)^{ms}  \sup_{z\in B_{\ell_q^n}}\sum_{\alpha \in\Lambda(m,n)}\frac{m!}{\alpha!} |z_1^s|^{\alpha_1}\cdots |z_n^s|^{\alpha_n} \\[2mm]
		&=&\norm{v}^s \norm{P}_{\mathcal{P}(^m \ell_p^n, X)}^s (en)^{ms}  \sup_{z\in B_{\ell_q^n}}\left(|z_1^s|+\cdots+ |z_n^s|\right)^m \\[2mm]
				&=&\norm{v}^s \norm{P}_{\mathcal{P}(^m \ell_p^n, X)}^s (en)^{ms}  \sup_{z\in B_{\ell_q^n}} \norm{z}_s^{ms} \\[2mm]
		&\leq& \norm{v}^s \norm{P}_{\mathcal{P}(^m \ell_p^n, X)}^s (en)^{ms}  \sup_{z\in B_{\ell_q^n}}\left(n^{\max\left\{0,  \frac{1}{s}-\frac{1}{q}\right\}}\norm{z}_q\right)^{ms} \\[2mm]
				&\leq& \lambda  e^{ms}  n^{ms\left(1+\max\left\{0,  \frac{1}{s}-\frac{1}{q}\right\}\right)}\norm{P}_{\mathcal{P}(^m \ell_p^n, X)}^s.
		\eeas
		
In view of \eqref{Defi3}, it follows that 
\beas
	K^s_m(B_{\ell^n_p}, B_{\ell^n_q},v,\lambda)\geq \frac{1}{e\, n^{\max\left\{1, 1+\frac{1}{s}-\frac{1}{q}\right\}}}\quad\mbox{for all $m\in\Bbb{N}$}.
\eeas
Using Lemma \ref{lem1} and for each $1\leq s<\infty$, we have 
\beas
K^s(B_{\ell^n_p},B_{\ell^n_q}, v, \lambda)\geq \left(\frac{\lambda-\norm{v}^s}{2\lambda-\norm{v}^s}\right)^{1/s}\frac{1}{e\, n^{\max\left\{1, 1+\frac{1}{s}-\frac{1}{q}\right\}}}.
\eeas
On the other hand, for every $m$-homogeneous polynomial $P(z)=\sum_{|\alpha|=m}c_{\alpha}z^{\alpha}$, we obtain 
\beas
	\sup_{z\in \frac{1}{\xi \sqrt[m]{\norm{v}}} B_{\ell^n_q}}\ \sum_{\alpha \in\Lambda(m,n)}\norm{v(c_{\alpha})z^{\alpha}}^s&=&
	\frac{1}{\norm{v}^s} \sup_{z\in \frac{1}{\xi}B_{\ell^n_q} }  \sum_{\alpha \in\Lambda(m,n)}\norm{v(c_{\alpha})z^{\alpha}}^s\\[2mm]&\leq &\norm{P}^s_{\mathcal{P}(^m\ell^n_p, X)},
\eeas
where $\xi=e\,n^{\max\left\{1, 1+\frac{1}{s}-\frac{1}{q}\right\}}$.
Thus, for every $m\in \Bbb{N}$, we have 
\beas
	K^s_m(B_{\ell^n_p},B_{\ell^n_q}, v)\geq \frac{1}{e\,n^{\max\left\{1, 1+\frac{1}{s}-\frac{1}{q}\right\}}\sqrt[m]{\norm{v}}}.
\eeas
Using Lemma \ref{lem1}, we obtain 
\beas
K^s(B_{\ell^n_p},B_{\ell^n_q}, v,\lambda)\geq 
\left\{\begin{array}{lll}
	&\left(\dfrac{\lambda-\norm{v}^s}{\lambda-\norm{v}^s+1}\right)^{1/s}\dfrac{1}{e\,n^{\max\left\{1, 1+\frac{1}{s}-\frac{1}{q}\right\}}\norm{v}}, &\mbox{if}\,\, \norm{v}\geq 1\\[5mm]
	&\left(\dfrac{\lambda-\norm{v}^s}{\lambda-\norm{v}^s+1}\right)^{1/s}\dfrac{1}{e\,n^{\max\left\{1, 1+\frac{1}{s}-\frac{1}{q}\right\}}}, &\mbox{if}\,\, 0<\norm{v}<1.
\end{array}\right.
\eeas
This completes the proof.
\end{proof}
As a consequence, we obtain the lower bound for the mixed Bohr radius of Banach spaces.
\begin{cor}\label{Cor1}
	Let $\lambda>1$, $1\leq p,q\leq \infty$ and $X$ be a complex Banach space. For $1\leq s<\infty$ and $n\in \Bbb{N}$, we have
	\beas
		K^s(B_{\ell^n_p}, B_{\ell^n_q},X,\lambda)\geq \left(\frac{\lambda-1}{\lambda}\right)^{1/s}\dfrac{1}{e\,n^{\max\left\{1, 1+\frac{1}{s}-\frac{1}{q}\right\}}}.
	\eeas
\end{cor}
\begin{proof} 
By substituting $\norm{v}=1$ in \textrm{Theorem} \ref{T1}, we have
\beas K^s(B_{\ell^n_p}, B_{\ell^n_q},v,\lambda)\geq \max \left\{\left(\dfrac{\lambda-1}{2\lambda-1}\right)^{1/s}, \left(\dfrac{\lambda-1}{\lambda}\right)^{1/s}\right\} \frac{1}{e\, n^{\max\left\{1, 1+\frac{1}{s}-\frac{1}{q}\right\}}}.\eeas
Using the fact that, $x\longrightarrow x^{1/s}$ is strictly increasing on $(0, \infty)$, we obtain desired conclusion. This completes the proof.\end{proof}
\noindent In the following result, we provide the asymptotic estimates for $K^s(B_{\ell^n_p},B_{\ell^n_q},X,\lambda)$, where $X$ is an infinite-dimensional Banach space. 
\begin{thm}\label{T2}
	Let $1\leq p,q\leq \infty$, $1\leq s<\infty$, $\lambda>1$ and $X$ be an infinite-dimensional complex Banach space of cotype $r$ with $2\leq r\leq \infty$. Then, we have 
	\beas
	\left(\frac{\lambda-1}{\lambda}\right)^{1/s}\frac{1}{e\,C_r(X)n^{\max\left\{\frac{1}{p},\,  1-\frac{1}{r}-\frac{1}{q}+\frac{1}{p}\right\}}}\leq K^s\left(B_{\ell^n_p},B_{\ell^n_q},X,\lambda\right)\leq \frac{\lambda^{1/s}}{n^{\min\left\{\frac{1}{s}-\frac{1}{q}, \frac{1}{s}-\frac{1}{q}-\frac{1}{\text{Cot}(X)}+\frac{1}{p}\right\}}}.\eeas
	\noindent Moreover, if $X$ has no finite cotype, then 
	\beas
		\left(\frac{\lambda-1}{\lambda}\right)^{1/s}\dfrac{1}{e\,n^{\max\left\{1, 1+\frac{1}{s}-\frac{1}{q}\right\}}}\leq K^s\left(B_{\ell^n_p},B_{\ell^n_q},X,\lambda\right)\leq \frac{\lambda^{1/s}}{n^{\min\left\{\frac{1}{s}-\frac{1}{q}, \frac{1}{s}-\frac{1}{q}+\frac{1}{p}\right\}}}.
	\eeas
\end{thm}
\begin{proof}
	In view of Maurey and Pisier \cite[Theorem 14.5]{DJT1995}, for every infinite-dimensional Banach space $X$, we have
	\beas
		Cot(X)=\sup \left\{2\leq r\leq \infty : X \quad\mbox{finitely factors}\quad\ell_r \hookrightarrow \ell_{\infty}\right\}.
	\eeas
Furthermore, the Banach space $X$ finitely factors $\ell_r \hookrightarrow \ell_\infty$ if for every $0<\delta<1$ and $n\in \Bbb{N}$, there are $x_1, \ldots, x_n\in X$ such that 
	\bea\label{a8} (1-\delta)\norm{z}_{\infty}\leq \norm{\sum_{k=1}^{n}z_k x_k}\leq \norm{z}_r\quad \mbox{for all $z=(z_1, \ldots, z_n)\in \ell_r^n$}.\eea
 In particular, for $z=e_k=(0,\ldots, 0,1,0,\ldots,0)$, we have
\beas
	1-\delta=(1-\delta)\norm{e_k}_{\infty}\leq \norm{x_k}.
\eeas
Hence, for $0<\delta<1$, there exist $x_1, \ldots, x_n\in X$ such that 
\bea\label{a9}
	\sum_{k=1}^{n}\norm{x_k}^s\geq n(1-\delta)^s,\quad\text{where}\quad 1\leq s<\infty.
\eea
For $x_1, \ldots, x_n\in X$, we define a linear polynomial $P(z)=\sum_{k=1}^n x_k z_k$. By Definition \ref{Defi2}, we have 
\bea\label{a10} \sup_{z\in B_{\ell^n_{q}}}\sum_{k=1}^n\norm{x_k z_k}^s\leq \frac{1}{\left(K^s_1\left(B_{\ell^n_p},B_{\ell^n_q},X\right)\right)^s}\sup_{z\in B_{\ell^n_{p}}}\norm{\sum_{k=1}^n x_k z_k}^s.\eea
Let $0<\epsilon<1$ and choose $z\in\Bbb{C}^n$ in such a way that $z_k=(1-\epsilon)/n^{1/q}$ for $1\leq k\leq n$. It is evident that $z\in B_{\ell^n_{q}}$. Applying (\ref{a10}) to this 
particular $z$, we obtain
\bea\label{a11} \frac{(1-\epsilon)^s}{n^{s/q}}\sum_{k=1}^n\norm{x_k}^s\leq \frac{1}{\left(K^s_1\left(B_{\ell^n_p},B_{\ell^n_q},X\right)\right)^s}\sup_{z\in B_{\ell^n_{p}}}\norm{\sum_{k=1}^n x_k z_k}^s.\eea
Letting $\epsilon\to 0^{+}$ in (\ref{a11}), and from (\ref{a8}) and (\ref{a9}), we obtain 
\beas
	\frac{n(1-\delta)^s}{n^{s/q}}\leq \frac{1}{n^{s/q}}\sum_{k=1}^n\norm{x_k}^s&\leq& \frac{1}{\left(K^s_1\left(B_{\ell^n_p},B_{\ell^n_q},X\right)\right)^s}\sup_{z\in B_{\ell^n_{p}}}\norm{\sum_{k=1}^n x_k z_k}^s\nonumber\\[2mm]
	&\leq& \frac{1}{\left(K^s_1\left(B_{\ell^n_p},B_{\ell^n_q},X\right)\right)^s}\sup_{z\in B_{\ell^n_{p}}}\norm{z}^s_{r}\nonumber\\[2mm]
	&\leq& \frac{1}{\left(K^s_1\left(B_{\ell^n_p},B_{\ell^n_q},X\right)\right)^s}\sup_{z\in B_{\ell^n_{p}}}\left(n^{\max\left\{0,  \frac{1}{r}-\frac{1}{p}\right\}}\norm{z}_p\right)^s\nonumber\\[2mm]
	&\leq &\left(\frac{n^{\max\left\{0,  \frac{1}{r}-\frac{1}{p}\right\}}}{K^s_1\left(B_{\ell^n_p},B_{\ell^n_q},X\right)}\right)^s.
\eeas 
Thus, for any $1\leq s<\infty$, we obtain 
\beas
	K^s_1\left(B_{\ell^n_p},B_{\ell^n_q},X\right) \leq \frac{n^{\max\left\{\frac{1}{q},  \frac{1}{q}+\frac{1}{r}-\frac{1}{p}\right\}}}{n^{1/s}(1-\delta)}
\eeas
holds for all $0<\delta<1$. Consequently, it follows that
\beas
	K^s_1\left(B_{\ell^n_p},B_{\ell^n_q},X\right)\leq \frac{1}{n^{\min\left\{\frac{1}{s}-\frac{1}{q}, \frac{1}{s}-\frac{1}{q}-\frac{1}{\text{Cot}(X)}+\frac{1}{p}\right\}}}.
\eeas
Using the fact that $K^s\left(B_{\ell^n_p},B_{\ell^n_q},X,\lambda\right)\leq K^s_1\left(B_{\ell^n_p},B_{\ell^n_q},X,\lambda\right)=\lambda^{1/s}K^s_{1}\left(B_{\ell^n_p},B_{\ell^n_q},X\right)$, we obtain the desired upper bound for $K^s\left(B_{\ell^n_p},B_{\ell^n_q},X,\lambda\right)$.\\ [2mm]
\indent Now we establish the lower bound for $K^s\left(B_{\ell^n_p},B_{\ell^n_q},X,\lambda\right)$, where $X$ is a Banach space of cotype $r$ with $2\leq r<\infty$. 
Let $P_m:\ell^n_p\rightarrow X$ be an $m$-homogeneous polynomial  such that $P_m(z)=\sum_{\alpha\in\Lambda(m,n)}c_{\alpha}z^{\alpha}$. Then there exists a unique 
symmetric $m$-linear mapping $A: \underbrace{\ell^n_p\times \cdots\times \ell^n_p}_{m\text{-times}}\rightarrow X$ such that $P_m(z)=A(z,\ldots,z)$ for every $z\in \Bbb{C}^n$. In view of \textrm{Lemma} \ref{lem3}, we have
\beas
	\left(\sum_{i\in \mathcal{M}(m,n)}\norm{A(e_{i_1}, \dots, e_{i_m})}^r\right)^{1/r}\leq C_r(X)^m n^{m/p}\norm{A},\eeas
	where $C_r(X)$ is the cotype $r$ constant of $X$. \\[2mm]
	Let $t>1$ be such that $1/r+1/t=1$.
Using the fact that $\norm{a}_s\leq \norm{a}_1$ for every $s\geq 1$ and every $a\in \Bbb{C}^n$, and applying the H\"{o}lder inequality, for all $z\in \Bbb{C}^n$, we obtain
\beas
	\sum_{\alpha \in\Lambda(m,n)}\norm{c_{\alpha}z^{\alpha}}^s
	& =& \sum_{\mathbf{j}\in \mathcal{J}(m,n)} \text{card}[\mathbf{j}]^s\norm{A(e_{i_1}, \ldots, e_{i_m})}^s |z_{i_1}\cdots z_{i_m}|^s\\[2mm]
	& \leq & \left(\sum_{\mathbf{j}\in \mathcal{J}(m,n)} \text{card}[\mathbf{j}]\norm{A(e_{i_1}, \ldots, e_{i_m})} |z_{i_1}\cdots z_{i_m}|\right)^s\\[2mm]
	&=& \left(	\sum_{i\in \mathcal{M}(m,n)} \norm{A(e_{i_1}, \ldots, e_{i_m})} |z_{i_1}\cdots z_{i_m}|\right)^s\\[2mm]
	&\leq &\left(\left(\sum_{i\in \mathcal{M}(m,n)}\norm{A(e_{i_1}, \ldots, e_{i_m})}^r\right)^{1/r}\left(\sum_{i\in \mathcal{M}(m,n)}|z_{i_1}\cdots z_{i_m}|^{t}\right)^{1/t}\right)^s\\[2mm]
	&\leq& C_r(X)^{ms} n^{ms/p} \norm{A}^s \left(\sum_{i\in \mathcal{M}(m,n)}|z_{i_1}\cdots z_{i_m}|^{t}\right)^{s/t}\\[2mm]
	&=&C_r(X)^{ms} n^{ms/p}\norm{A}^s \left(|z_1|^t+\cdots+|z_n|^t\right)^{\frac{ms}{t}}\\[2mm]
		&\leq &C_r(X)^{ms} n^{ms/p}\norm{A}^s n^{ms \max\left\{0,\,  1-\frac{1}{r}-\frac{1}{q}\right\}}\norm{z}_q^{ms}.
\eeas
Using a consequence of the polarization formula, we obtain 
\beas
\sup_{z\in B_{\ell_q^n}}\sum_{\alpha \in\Lambda(m,n)}\norm{c_{\alpha}z^{\alpha}}^s &\leq& C_r(X)^{ms} n^{ms/p} \norm{A}^s n^{ms \max\left\{0, \, 1-\frac{1}{r}-\frac{1}{q}\right\}}\\[2mm]
&\leq&  C_r(X)^{ms} n^{ms \left(\frac{1}{p}+\max\left\{0, \, 1-\frac{1}{r}-\frac{1}{q}\right\}\right)}\left(\frac{m^m}{m!}\right)^s \norm{P_m}^s_{B_{\ell^n_p}} \\[2mm]
&\leq&  C_r(X)^{ms}n^{ms \max\left\{\frac{1}{p}, \, 1-\frac{1}{r}-\frac{1}{q}+\frac{1}{p}\right\}}e^{ms} \norm{P_m}^s_{B_{\ell^n_p}}.
\eeas
Thus, for every $m\in \Bbb{N}$ and $1\leq s<\infty$, we obtain
\beas
K^s_m\left(B_{\ell^n_p}, B_{\ell^n_q},X\right) \geq \frac{1}{e\,C_r(X)n^{\max\left\{\frac{1}{p}, \, 1-\frac{1}{r}-\frac{1}{q}+\frac{1}{p}\right\}}}.
\eeas
Using \textrm{Lemma} \ref{lem1}, we obtain
 \beas
 K^s\left(B_{\ell^n_p}, B_{\ell^n_q},X,\lambda\right) \geq \left(\frac{\lambda-1}{\lambda}\right)^{1/s}\frac{1}{e\,C_r(X)n^{\max\left\{\frac{1}{p},\,  1-\frac{1}{r}-\frac{1}{q}+\frac{1}{p}\right\}}}\eeas
for every $\lambda>1$, $1\leq p,q\leq \infty$, and $1\leq s< \infty$. 
 For Banach spaces $X$ without finite cotype, the lower bound for $K^s\left(B_{\ell^n_p}, B_{\ell^n_q},X,\lambda\right)$ has been established in \textrm{Corollary} \ref{Cor1}.
This completes the proof.
\end{proof}
As an application of \textrm{Theorem} \ref{T2}, we obtain the following examples.
\begin{ex}
Let $X = Y = \ell_2$ and let $v = \textrm{id}:\ell_2\to \ell_2$. Then, we have $\text{Cot}(\ell_2) = 2$.
For $p=q = r = 2$,  $s= 1$ and $\lambda = 2$,  \textrm{Theorem} \ref{T2} gives
\beas
\frac{1}{2e\,C_2(\ell_2)\,n^{1/2}}\leq K^1(B_{\ell_2^n}, B_{\ell_2^n}, \ell_2, 2)\leq \frac{2}{n^{1/2}}.
\eeas
Hence, we obtain the sharp asymptotic order $K^1(B_{\ell_2^n}, B_{\ell_2^n}, \ell_2, 2) \sim n^{-1/2}$.
\end{ex}
\begin{ex}
Let $X = Y = \ell_2$ and let $v = \textrm{id}:\ell_2\to \ell_2$. Then, we have $\text{Cot}(\ell_2) = 2$.
For $p=2$, $q =3$, $r=2$, $s= 2$ and $\lambda = 2$,  \textrm{Theorem} \ref{T2} gives
\beas
\frac{1}{\sqrt{2}\,e\,C_2(\ell_2)\,n^{2/3}}\leq K^2(B_{\ell_2^n}, B_{\ell_3^n}, \ell_2, 2)\leq\frac{\sqrt{2}}{n^{1/6}}.
\eeas
Since $n^{-2/3} \ll n^{-1/6}$ as $n \to \infty$, the lower and upper bounds do not match, and the sharp asymptotic order of 
$K^2(B_{\ell_2^n}, B_{\ell_3^n}, \ell_2, 2)$ is not determined by \textrm{Theorem} \ref{T2} 
for this combination of parameters. The actual decay rate lies somewhere between $n^{-2/3}$ and $n^{-1/6}$.
\end{ex}
\begin{ex}
Let $X = Y = \ell_2$ and let $v = \textrm{id}:\ell_2\to \ell_2$. Then, we have $\text{Cot}(\ell_2) = 2$.
For $p = r = 2$,  $q = 4$, $\lambda = 2$ and $s = 2$,  \textrm{Theorem} \ref{T2} gives
\beas
\frac{1}{\sqrt{2}\,e\,C_2(\ell_2)\; n^{3/4}}\leq K^2(B_{\ell_2^n}, B_{\ell_4^n}, \ell_2, 2)\leq\frac{\sqrt{2}}{n^{1/4}}.
\eeas
\end{ex}
\vspace{2mm}
\noindent In $2012$, Defant {\it et al.} \cite{DMS2012} established the following result, which is one of the main tools used to estimate the multidimensional Bohr radii of operators. 
 \begin{theoB} \cite[Theorem 5.3]{DMS2012} Let $Y$ be an $r$-concave Banach lattice with $2\leq r<\infty$, and $v:X\rightarrow Y$ an $(t,1)$-summing operator with $1\leq t\leq r$. Then there is a constant $C>0$ such that for every $m$-homogeneous polynomial $P:\ell^n_{\infty}\rightarrow X$ with $P(z)=\sum_{\alpha \in\Lambda(m,n)}c_{\alpha}z^{\alpha}$, we have  
\beas
	\left(\sum_{\alpha \in\Lambda(m,n)}\norm{v(c_{\alpha})}^{\rho}_{Y}\right)^{1/\rho}\leq C^m\norm{P}_{\mathcal{P}(^m \ell^n_{\infty},X)},\quad\mbox{where $\rho=\dfrac{rtm}{r+(m-1)t}$}.
\eeas
\end{theoB}
\noindent Using the diagonal operator $D_{z}:\ell^n_{\infty}\rightarrow \ell^n_{p}$ ($z\in B_{\ell_p^n}$), Kumar and Manna \cite{KM2025} have established the following result for $m$-homogeneous polynomials $P\in \mathcal{P}(^m\ell^n_p,X)$.
\begin{theoC}\cite{KM2025}\label{TheoC}
Suppose $Y$ is an $r$-concave Banach lattice with $2\leq r<\infty$ or Banach space of cotype $r$ with $2\leq r<\infty$ and $v:X\rightarrow Y$ is an $(t,1)$-summing operator
with $1\leq t\leq r$. Then there is a constant $C>0$ such that for every $m$-homogeneous polynomial $P(z)=\sum_{\alpha \in\Lambda(m,n)}c_{\alpha}z^{\alpha}\in \mathcal{P}(^m\ell^n_{p}, X)$, we have  
	\beas
		\left(\sum_{\alpha \in\Lambda(m,n)}\norm{v(c_{\alpha})z^{\alpha}}^{\rho}_{Y}\right)^{1/\rho}\leq {C}^m\norm{P}_{\mathcal{P}(^m \ell^n_{p}, X)},\quad\mbox{where $\rho=\dfrac{rtm}{r+(m-1)t}$}.
	\eeas
\end{theoC}
\noindent In the following result, we derive the lower estimates for $K^s(B_{\ell^n_p}, B_{\ell^n_q},v,\lambda)$ for any bounded linear operator $v$.
\begin{thm}\label{T3}
	Let $1\leq p,q\leq \infty$, $1\leq s<\infty$ and $v:X\rightarrow Y$ be a bounded operator between Banach spaces $X$ and $Y$ with $\norm{v}<\lambda^{1/s}$.
	\begin{enumerate}
		\item[(i)] Let $X$ or $Y$ be of cotype $r$ with $2\leq r\leq \infty$. Then there is a constant $C>0$ such that for every $\norm{v}<\lambda^{1/s}$ and every $n$, we have 
		\beas
		K^s(B_{\ell^n_p}, B_{\ell^n_q},v,\lambda) \geq C \left(\frac{\lambda - \norm{v}^s}{\lambda}\right)^{1/s} \frac{1}{n^{\max\left\{\frac{1}{p},\,  1-\frac{1}{r}-\frac{1}{q}+\frac{1}{p}\right\}}}.
		\eeas
		\item[(ii)] Let $Y$ be an $r$-concave Banach lattice with $2\leq r<\infty$ and there is a $t\in[1,r)$ such that $v$ is a $(t,1)$-summing operator. Then there is a constant $D>0$  such that for every $\norm{v}<\lambda^{1/s}$ and every $n$, we have 
		\beas
			K^s(B_{\ell^n_p}, B_{\ell^n_q},v,\lambda)  \geq D\, \left(\frac{\lambda-\norm{v}^s}{2\lambda -\norm{v}^s}\right)^{1/s} \frac{1}{n^{\max\left\{0,\, \frac{1}{p}-\frac{1}{q}\right\}}} \left(\frac{\log n}{n}\right)^{1-\frac{1}{r}}.
		\eeas

	\end{enumerate}
\end{thm}
\begin{proof}
	(i) From (\ref{a3}) and \textrm{Theorem} \ref{T2}, we have 
	\beas K^s(B_{\ell^n_p},B_{\ell^n_q},v,\lambda)&\geq& \max \left\{K^s\left(B_{\ell^n_p},B_{\ell^n_q},X,\lambda/\norm{v}^s\right), K^s\left(B_{\ell^n_p},B_{\ell^n_q},Y,\lambda/\norm{v}^s\right)\right\}\\[2mm]
	&\geq& \left(\frac{\lambda-\norm{v}^s}{\lambda}\right)^{1/s}\frac{1}{e\,\min\{C_r(X), C_r(Y)\}n^{\max\left\{\frac{1}{p},\,  1-\frac{1}{r}-\frac{1}{q}+\frac{1}{p}\right\}}}.
	\eeas
Hence, the result follows.\\[1mm]
(ii)   Let $P(w)=\sum_{\alpha \in\Lambda(m,n)}c_{\alpha} w^{\alpha}\in \mathcal{P}(^m \ell^n_p,X)$ be an $m$-homogeneous polynomial and $z\in B_{\ell_q^n}$. Using the norm comparison, we have $\norm{z}_p\leq n^{\max\left\{0,\, \frac{1}{p}-\frac{1}{q}\right\}}\norm{z}_q< n^{\max\left\{0,\, \frac{1}{p}-\frac{1}{q}\right\}}$.
We can write $z$ as $z=\mu w$ with $\mu=n^{\max\left\{0,\, \frac{1}{p}-\frac{1}{q}\right\}}$ and $w=z/\mu \in B_{\ell_p^n}$.
Thus, we obtain
\beas
\sum_{\alpha\in\Lambda(m,n)}\norm{v(c_\alpha) z^\alpha}^s=\mu^{ms}\sum_{\alpha\in\Lambda(m,n)}\norm{v(c_\alpha) w^\alpha}^s&=& n^{ms\max\left\{0,\, \frac{1}{p}-\frac{1}{q}\right\}}\sum_{\alpha\in\Lambda(m,n)}\norm{v(c_\alpha) w^\alpha}^s\\[2mm]
&\leq& n^{ms\max\left\{0,\, \frac{1}{p}-\frac{1}{q}\right\}}\left(\sum_{\alpha\in\Lambda(m,n)}\norm{v(c_\alpha) w^\alpha}\right)^s . \eeas
Applying \textrm{Theorem C} and the H\"{o}lder inequality, there is a constant $C>0$ such that 
 \beas
 	\sum_{\alpha \in\Lambda(m,n)}\norm{v(c_{\alpha})z^{\alpha}}^s &\leq &n^{ms\max\left\{0,\, \frac{1}{p}-\frac{1}{q}\right\}}\left(\left(\sum_{\alpha \in\Lambda(m,n)}1\right)^{1-1/\rho}\left(\sum_{\alpha \in\Lambda(m,n)}\norm{v(c_{\alpha})w^{\alpha}}^{\rho}\right)^{1/\rho}\right)^s\\[2mm] 
 	&\leq&n^{ms\max\left\{0,\, \frac{1}{p}-\frac{1}{q}\right\}} \left(\sum_{\alpha \in\Lambda(m,n)}1\right)^{\frac{((r-1)mt-r+t)s}{rtm}} C^{ms} \norm{P}^s_{\mathcal{P}(^m \ell^n_{p},X)},
 \eeas
where $\rho=rtm/(r+(m-1)t)$. Note that 
\bs\beas e^m\geq \frac{m^m}{m!}\geq  \frac{m^m}{m!}\cdot\frac{n(n+1)\cdots(n+m-1)}{(n+m)^m}=\frac{m^m}{(n+m)^m}\cdot\frac{(n+m-1)!}{(n-1)! m!}=\binom{n+m-1}{m}\left(\frac{m}{n+m}\right)^m.\eeas\es
Thus, we have 
\beas
\sum_{\alpha \in\Lambda(m,n)}1=\binom{n+m-1}{m}\leq e^m\left(1+\frac{n}{m}\right)^m
\eeas
and consequently,
\beas
 	\sup_{z\in B_{\ell_q^n}}\sum_{\alpha \in\Lambda(m,n)}\norm{v(c_{\alpha})z^{\alpha}}^s \leq n^{ms\max\left\{0,\, \frac{1}{p}-\frac{1}{q}\right\}}\left(e\left(1+\frac{n}{m}\right)\right)^{\frac{((r-1)mt-r+t)ms}{mrt}} C^{ms} \norm{P}^s_{\mathcal{P}(^m \ell^n_{p},X)}.\eeas
In view of \eqref{Defi3} and \eqref{a2}, there is a constant $D_1>0$ such that 
\bea\label{a12}
	K^s_m\left(B_{\ell^n_p},B_{\ell^n_q},v,\lambda\right)&\geq& \frac{D_1\, \lambda^{1/(ms)}}{n^{\max\left\{0,\, \frac{1}{p}-\frac{1}{q}\right\}}}\left(1+\frac{n}{m}\right)^{-\frac{(r-1)mt-r+t}{rmt}}\nonumber\\[2mm]
	&\geq& \frac{D_1\, \min\{1,\lambda^{1/s}\}}{n^{\max\left\{0,\, \frac{1}{p}-\frac{1}{q}\right\}}}\left(1+\frac{n}{m}\right)^{-\frac{(r-1)mt-r+t}{rmt}}\eea
for every $m$, $1\leq s<\infty$, and $\norm{v}<\lambda^{1/s}$. 
Since $t\in[1,r)$, we observe that 
\beas \left(1+\frac{n}{m}\right)^{\frac{(r-1)mt-r+t}{rmt}}=\left(1+\frac{n}{m}\right)^{\left(1-\frac{1}{r}\right)-\frac{1}{m}\left(\frac{1}{t}-\frac{1}{r}\right)}= \left(1+\frac{n}{m}\right)^{\beta_1-\frac{\beta_2}{m}},\eeas
where $\beta_1=1-1/r> 0$ and $\beta_2=1/t-1/r>0$. Thus, we have 
\beas \left(1+\frac{n}{m}\right)^{-\left(\beta_1-\frac{\beta_2}{m}\right)}&\geq& 2^{-\left(\beta_1-\frac{\beta_2}{m}\right)}\min\left\{1, \left(\frac{n}{m}\right)^{-\left(\beta_1-\frac{\beta_2}{m}\right)}\right\}\\[2mm]
&\geq&\frac{1}{2^{\beta_1}}\min\left\{1, \left(\frac{m n^{\frac{\beta_2}{\beta_1 m}}}{n}\right)^{\beta_1}\cdot \frac{1}{m^{\beta_2/m}}\right\}.\eeas
Let us consider the function $g: (0, \infty)\to \Bbb{R}$ defined by $g(x)= x n^{\beta_2/(\beta_1 x)}$, where $n\geq 2$. It is evident that 
\beas g'(x)=n^{\frac{\beta_2}{\beta_1 x}}\left(1-\frac{\beta_2 \log n}{\beta_1 x}\right).\eeas
It follows that the function $g$ is strictly decreasing on $(0, \beta_2 \log n/\beta_1)$ and strictly increasing on $(\beta_2 \log n/\beta_1, \infty)$. Therefore, the function $g$ attains its unique global minimum at $x=\beta_2 \log n/\beta_1$ and consequently, for every $m\in\Bbb{N}$, we have 
\beas m n^{\frac{\beta_2}{\beta_1 m}}\geq \min_{x>0} g(x)=\frac{e\beta_2 \log n}{\beta_1}.\eeas
Similarly, the function $h: (0, \infty)\to \Bbb{R}$ defined by $h(x)= x^{-\beta_2/x}$ is strictly decreasing on $(0, e)$ and strictly increasing on $(e, \infty)$. Therefore, the function 
$h$ attains its unique global minimum at $x=e$ and consequently, for every $m\in\Bbb{N}$, we have $m^{-\beta_2/m}\geq \min_{x>0} h(x)=e^{-\beta_2/e}$. Thus, we have
\beas \left(1+\frac{n}{m}\right)^{-\left(\beta_1-\frac{\beta_2}{m}\right)}\geq
\frac{1}{2^{\beta_1}}\min\left\{1, \left(\frac{e\log n}{ n}\right)^{\beta_1}\cdot \frac{\beta_2^{\beta_1} }{e^{\frac{\beta_2}{e}} \beta_1^{\beta_1}}\right\}.\eeas
Using the fact that $e^{\frac{\beta_2}{e}}>1$, $\beta_2\leq \beta_1$ and $(e\log n)/n<1$ for all $n\in\Bbb{N}$, we have 
\bea\label{a13} \left(1+\frac{n}{m}\right)^{-\left(\beta_1-\frac{\beta_2}{m}\right)}\geq \frac{1}{2^{\beta_1}}\left(\frac{e\log n}{ n}\right)^{\beta_1}\cdot \frac{\beta_2^{\beta_1} }{e^{\frac{\beta_2}{e}} \beta_1^{\beta_1}}=D_2\left(\frac{\log n}{ n}\right)^{\beta_1},\eea
which holds for all $m,n$, where $D_2=(e\beta_2)^{\beta_1}/\left(2^{\beta_1}e^{\frac{\beta_2}{e}} \beta_1^{\beta_1}\right)$. For $n=1$ and $m\in\Bbb{N}$, we have 
\bea\label{a14}\left(1+\frac{n}{m}\right)^{-\left(\beta_1-\frac{\beta_2}{m}\right)}\geq\frac{1}{2^{\beta_1}}\min\left\{1, \frac{1}{m^{\beta_2/m}}\right\}\geq \frac{1}{2^{\beta_1}}\min\left\{1, \frac{1}{e^{\beta_2/e}}\right\}=\frac{1}{2^{\beta_1}e^{\beta_2/e}}.\eea

From (\ref{a12}), (\ref{a13}) and (\ref{a14}), we see that there is a constant $D>0$ such that for every $n\in\Bbb{N}$, $1\leq s<\infty$ and $\norm{v}<\lambda^{1/s}$, we obtain
\beas
	\inf\limits_{m\in \Bbb{N}} \left\{K^s_m\left(B_{\ell^n_p},B_{\ell^n_q},v,\lambda\right)\right\}\geq \frac{D}{n^{\max\left\{0,\, \frac{1}{p}-\frac{1}{q}\right\}}} \left(\frac{\log n}{n}\right)^{1-\frac{1}{r}}.\eeas
Using \textrm{Lemma} \ref{lem1}, we obtain desired conclusion. This completes the proof.\end{proof}
\begin{rem}
In the special case $r = 2$ and $p \geq q$ with $p \geq 2$ of \textrm{Theorem \ref{T3}} (ii), then the lower bound becomes
\beas
K^s(B_{\ell_p^n}, B_{\ell_q^n}, v, \lambda)\ge D \left(\frac{\lambda - \|v\|^s}{2\lambda - \|v\|^s}\right)^{1/s}\left(\frac{\log n}{n}\right)^{1/2}.
\eeas
The exponent $1/2$ coincides with the sharp exponent $1 - \frac{1}{\min\{p,2\}} = \frac{1}{2}$ in the Defant--Frerick asymptotics for the scalar Bohr radius $K(B_{\ell_p^n})$
(see \cite[Theorem 1.1]{DF2011}).
\end{rem}
\section{Declarations}
\noindent{\bf Acknowledgment:} The authors would like to express their sincerest gratitude to the referee for careful reading of the manuscript and many valuable suggestions, 
which greatly helped to improve the clarity of the exposition in this manuscript. The research of the first named author is supported by NBHM Research Grant (No. 02011/19/2025/NBHM (R.P.)/R\&D II/9461), DAE, Govt. of India. \\[1mm]
{\bf Conflict of Interest:} The authors declare that there are no conflicts of interest regarding the publication of this paper.\\[1mm]
{\bf Availability of data and materials:} Not applicable.\\[1mm]
{\bf Authors' contributions:} All authors contributed equally to the investigation of the problem, and all authors have read and approved the final manuscript.


\begin{thebibliography}{99}

	\bibitem{AAP2017}{\sc Y. Abu-Muhanna, R. M. Ali} and {\sc S. Ponnusamy}, On the Bohr inequality, In: N.K. Govil {\it et al.} (eds), 
Progress in Approximation Theory and Applicable Complex Analysis, {\it Springer Optimization and Its Applications} {\bf 117} (2017), 265--295.
	\bibitem{A2000} {\sc L. Aizenberg}, Multidimensional analogues of Bohr's theorem on power series, \textit{Proc. Amer. Math. Soc.} {\bf 128} (2000), 1147--1155.
	
	\bibitem{AAD2000} {\sc L. Aizenberg, A. Aytuna} and {\sc P. Djakov}, An abstract approach to Bohr's phenomenon, {\it Proc. Amer. Math. Soc.} {\bf 128} (2000), 2611--2619.
	
	
	\bibitem{AAD2001} {\sc L. Aizenberg, A. Aytuna}  and {\sc P. Djakov}, Generalization of theorem on Bohr for bases in spaces of holomorphic functions of several complex variables, {\it J. Math. Anal. Appl.} {\bf  258} (2001), 429--447.
	
	

	
\bibitem{ABM2026} {\sc V. Allu, R. Biswas} and {\sc R. Mandal}, The Bohr's phenomenon involving multiple Schwarz functions, {\it J. Math. Anal. Appl.} {\bf553} (2026), 129965.
	
	\bibitem{AHP2026} {\sc V. Allu, H. Halder} and {\sc S. Pal}, Multidimensional Bohr Radii for Holomorphic Functions with Values in Complex Banach Spaces, {\it Complex Anal. Oper. Theory} {\bf20} (2026), 97. 
	
	\bibitem{AP2026} {\sc V. Allu} and {\sc S. Pal}, On multidimensional Bohr radii for Banach spaces, {\it J. Oper. Theory} {\bf95} (2) (2026), 481--503.

	
\bibitem{BCQ2006} {\sc R. Balasubramanian, B. Calado} and {\sc H. Queff\'{e}lec}, The Bohr inequality for ordinary Dirichlet series, {\it Studia Math.} {\bf 175} (2006), 285--304. 

	
	\bibitem{BPS2014} {\sc F. Bayart, D. Pellegrino} and {\sc J. B. Seoane-Sep\'ulveda}, The Bohr radius of the $n$-dimensional polydisk is equivalent to $\sqrt{(log \, n)/n}$, {\it Adv. Math.} {\bf 264} (2014), 726--746.
	
	\bibitem{BDS2019} {\sc F. Bayart, A. Defant} and {\sc S. Schl\"{u}ters}, Monomial convergence for holomorphic functions on $\ell_r$, {\it J. Anal. Math.} {\bf 138} (2019), 107--134.
	

	\bibitem{BDK2004} {\sc C. B{\'e}n{\'e}teau}, {\sc A. Dahlner} and {\sc D. Khavinson}, Remarks on the Bohr phenomenon, {\it Comput. Methods Funct. Theory} \textbf{4}(1) (2004), 1--19.
	

	
	\bibitem{B2010} {\sc O. Blasco}, The Bohr radius of a Banach space, In G. P. Curbera {\it et al.} (eds), Vector measures, integration and related topics, {\it Oper. Theory Adv. Appl.}, {\bf 201}, Birkh{\"a}user Basel, 2010.
	
	\bibitem{B2017} {\sc O. Blasco}, The $p$-Bohr radius of a Banach space, {\it Collect. Math.} {\bf 68} (2017), 87--100.
	
	
	\bibitem{BK1997} {\sc H. P. Boas} and {\sc D. Khavinson}, Bohr's power series theorem in several variables, {\it Proc. Amer. Math. Soc.}  {\bf 125} (1997), 2975--2979.
	
		\bibitem{B2000} {\sc H. P. Boas}, Majorant Series,  {\it J. Korean Math. Soc.}  {\bf 37} (2000), 321--337.
	
	\bibitem{B1914} {\sc H. Bohr}, A theorem concerning power series,  {\it Proc. Lond. Math. Soc.} s2-13 (1914), 1--5.
	\bibitem{BH1931}{\sc H. F. Bohnenblust} and {\sc E. Hille}, On the absolute convergence of Dirichlet series, {\it Ann. of Math.} {\bf32}(3) (1931), 600--622.
	\bibitem{BPV2004} {\sc F. Bombal, D. P\'{e}rez-Garc\'{i}a} and {\sc I. Villanueva}, Multilinear extensions of Grothendieck's theorem, {\it Q. J. Math} {\bf 55} (2004), 441--450.
	
	\bibitem{B1962} {\sc E. Bombieri}, Sopra un teorema di H. Bohr e G. Ricci sulle funzioni maggioranti delle serie di potenze, {\it Boll. Un. Mat. Ital.} {\bf 17} (1962), 276--282.
	
	\bibitem{BB2004} {\sc E. Bombieri} and {\sc J. Bourgain}, A remark on Bohr's inequality, {\it Internat. Math. Res. Notices} {\bf 80} (2004), 4307--4330.
	\bibitem{DLP2025} {\sc H. G. Dai, M. S. Liu} and {\sc S. Ponnusamy}, The improved and refined Bohr type inequalities of the Dirichlet series, {\it Bull. Malays. Math. Sci. Soc.} {\bf 48} (2025), 79.

	
	\bibitem{D2023} {\sc N. Das}, Estimates for generalized Bohr radii in one and higher dimensions, {\it Canad. Math. Bull.} {\bf 66} (2023), 682--699.
	
	\bibitem{D2024} {\sc N. Das}, The $p$-Bohr radius for vector-valued holomorphic and pluriharmonic functions, {\it Forum Math.} {\bf 36} (2024), 765--782.
	
	\bibitem{DF2006} {\sc A. Defant} and {\sc L. Frerick}, A logarithmic lower bound for multi-dimensional Bohr radii, {\it Israel J. Math.} {\bf 152} (2006), 17--28.
	
   \bibitem{DF2011} {\sc A. Defant} and {\sc L. Frerick}, The Bohr radius of the unit ball of $\ell^n_p$, {\it J. Reine Angew. Math.} {\bf 660} (2011), 131--147.
	
	\bibitem{DFOOS2011} {\sc A. Defant, L. Frerick, J. Ortega-Cerd${\rm \grave{A}}$, M. Ouna{\"i}es} and {\sc K. Seip}, The Bohnenblust-Hille inequality for homogeneous polynomials is hypercontractive, {\it Ann. of Math.} {\bf 174} (2011), 512--517.
	
	\bibitem{DGM2003} {\sc A. Defant, D. Garc\'{i}a} and {\sc  M. Maestre}, Bohr power series theorem and local Banach space theory, {\it J. Reine Angew. Math.} {\bf 557} (2003), 173--197.

	
	\bibitem{DMP2009} {\sc A. Defant,  M. Maestre} and {\sc  C. Prengel}, Domains of convergence for monomial expansions of holomorphic functions in infinitely many variables, {\it J. Reine Angew. Math.} {\bf 634} (2009), 13--49.
	
	\bibitem{DMS2012} {\sc A. Defant,  M. Maestre} and {\sc U. Schwarting}, Bohr radii of vector valued holomorphic functions, {\it Adv. Math.} {\bf 231} (2012), 2837--2857.

\bibitem{DJT1995} {\sc J. Diestel, H. Jarchow} and {\sc A. Tonge}, Absolutely Summing Operators, Cambridge Studies in Advanced Mathematics, vol. {\bf 43}, Cambridge University Press, Cambridge, 1995.
	
	
	\bibitem{D1995} {\sc P. G. Dixon}, Banach algebras satisfying the non-unital von Neumann inequality, \textit{Bull. Lond. Math. Soc.} \textbf{27} (4) (1995), 359--362.
	
	\bibitem{DR2000} {\sc P. B. Djakov} and {\sc M. S. Ramanujan}, A remark on Bohr's theorem and its generalizations, \textit{J. Anal.} \textbf{8} (2000), 65--77.
	
	\bibitem{DT1989} {\sc S. Dineen} and {\sc R. M. Timoney}, Absolute bases, tensor products and a theorem of Bohr, \textit{Studia Math.} \textbf{94} (1989), 227--234.
	

	\bibitem{GMM2020} {\sc D. Galicer, M. Mansilla} and {\sc S. Muro}, Mixed Bohr radius in several variables, {\it Trans. Amer. Math. Soc.} {\bf 373} (2020), 777--796.
	
	\bibitem{HHK2009} {\sc H. Hamada, T. Honda} and {\sc G. Kohr}, Bohr's theorem for holomorphic mappings with values in homogeneous balls, {\it Israel J. Math.} {\bf 173} (2009), 177--187.

	\bibitem{KP2019}{\sc  I. R. Kayumov} and {\sc S. Ponnusamy}, On a powered Bohr inequality, {\it Ann. Acad. Sci. Fenn. Math.} {\bf44} (2019), 301--310.
	\bibitem{KM2025} {\sc S. Kumar} and {\sc R. Manna}, Multi-dimensional Bohr radii of Banach space valued holomorphic functions, {\it Bull. Sci. Math.} {\bf 205} (2025), 103688.
		\bibitem{KP2026} {\sc S. Kumar} and {\sc S. Ponnusamy}, Bounds on mixed Bohr radii of vector-valued holomorphic functions on Banach spaces, {\it Anal. Math. Phys.} {\bf16} (2026), 27.
		\bibitem{KPW2025} {\sc S. Kumar, S. Ponnusamy} and {\sc G. B. Williams},  The Bohr-type inequalities for holomorphic functions with lacunary series in complex Banach space,  {\it New York J. Math.} {\bf31} (2025), 259--281.
   \bibitem{LT1977} {\sc J. Lindenstrauss} and {\sc L. Tzafriri}, {\it Classical Banach Spaces I: Sequence Spaces}, in: Ergebnisse der Mathematik und ihrer Grenzgebiete, {\bf 92}, Springer-Verlag, Berlin-Heidelberg-New York, 1977.
   
    \bibitem{LT1979} {\sc J. Lindenstrauss} and {\sc L. Tzafriri}, {\it Classical Banach Spaces II: Function Spaces}, in: Ergebnisse der Mathematik und ihrer Grenzgebiete, {\bf 97}, Springer-Verlag, Berlin-Heidelberg-New York, 1979.
   \bibitem{LLP2021} {\sc G. Liu, Z. H. Liu} and {\sc S. Ponnusamy}, Refined Bohr inequality for bounded analytic functions, {\it Bull. Sci. Math.} {\bf173} (2021), 103054.
	
\bibitem{LP2021} {\sc M. S. Liu} and {\sc S. Ponnusamy}, Multidimensional analogues of refined Bohr's inequality, {\it Proc. Amer. Math. Soc.} {\bf 149} (2021), 2133--2146.
	

	
	\bibitem{PPS2002} {\sc V. I. Paulsen, G. Popescu} and {\sc D. Singh}, On Bohr's inequality, {\it Proc. Lond. Math. Soc.} s3-85 (2002), 493--512.
	
	

	\bibitem{P2019} {\sc G. Popescu}, Bohr inequalities for free holomorphic functions on polyballs, {\it Adv. Math.} {\bf 347} (2019), 1002--1053.
	
	
\end{thebibliography}
\end{document}